\documentclass[11pt]{article}%
\usepackage{amsmath}
\usepackage{amsfonts}
\usepackage{amssymb}
\usepackage{graphicx}
\usepackage{theorem}
\usepackage{color}%
\providecommand{\U}[1]{\protect\rule{.1in}{.1in}}
\newtheorem{thm}{Theorem}[section]
\newtheorem{lm}[thm]{Lemma}
\newtheorem{pr}[thm]{Proposition}
\newtheorem{df}[thm]{Definition}
\newtheorem{rmk}[thm]{Remark}
\newtheorem{cor}[thm]{Corollary}
{\theorembodyfont{\upshape}
\newtheorem{examp}[thm]{Example}
}
\numberwithin{equation}{section} 
\begin{document}

\title{Title: }

\begin{center}

\textbf{DUALLY CONE--BOUNDEDNESS OF\ A\ SET\ AND\ APPLICATIONS}

\bigskip

by

\bigskip

Marius\ DUREA\footnote{{\small Faculty of Mathematics, \textquotedblleft
Alexandru Ioan Cuza\textquotedblright\ University, 700506--Ia\c{s}i, Romania
and \textquotedblleft Octav Mayer\textquotedblright\ Institute of Mathematics,
Ia\c{s}i Branch of Romanian Academy, 700505--Ia\c{s}i, Romania{; e-mail:
\texttt{durea@uaic.ro}}}} and Elena-Cristina
STAMATE\footnote{{\small \textquotedblleft Octav Mayer\textquotedblright%
\ Institute of Mathematics, Ia\c{s}i Branch of Romanian Academy,
700505--Ia\c{s}i, Romania{; e-mail: }\texttt{cristina.stamate@acadiasi.ro}}}
\end{center}

\bigskip

\noindent{\small {\textbf{Abstract:}} We introduce and study a generalized
concept of boundedness of a subset of a normed vector space with respect to a
cone, which is defined as lower boundedness of the images of the underlying
set through all the positive functionals of the cone. We show that this is a
weaker notion when compared to other similar ones and we explore several links
with the existing literature. We subsequently demonstrate that this concept
furnishes the properties required to obtain various generalizations of
important results and techniques, including conic cancellation rules and the
R\aa dstr\"{o}m embedding procedure.}

\bigskip



\begin{center}

\end{center}

\section{Introduction and preliminaries}

Recently, conic notions and possible extensions of classical results to conic
settings have become a fruitful direction of research (see \cite{Bau},
\cite{DF-coneComp}, \cite{Far}, \cite{EE2} and the references therein). Of
course, the conic-related topological concepts and their potential
applications are not entirely new (see \cite{PT}, \cite{Luc}, for instance),
but the recent advances in the study of vector and set optimization problems
have given new impetus to further exploration of these ideas.

This paper primarily examines certain notions of boundedness of a set with
respect to a cone, and its objective is to relax the most common of these
notions in such a way that the main ingredients needed for applications in
optimization problems are still preserved by the weaker concept we introduce.
We identify such an concept and we are able to show that this is the weakest
one can provide for the validity of some important results as the conic
cancellation laws (see \cite{DFcancellation}). Surprisingly, this approach
also brings into focus classes of sets that have appeared in diverse contexts
in recent literature such as hyperbolic sets and pseudo-hyperbolic sets (see
\cite{Bair85}, \cite{Gos}, \cite{MP}). A careful comparison of these closely
related concepts is presented and their main differences are illustrated
through examples. Several features of the newly introduced notion allow us to
weaken the assumptions of recent results in literature (e.g., some results in
\cite{DSta}) concerning classes of optimization problems, including vector and
set optimization. These findings underscore the need to revisit existing
notions and motivate a deeper investigation of cone-bounded sets, both for
theoretical insight and applications in optimization.

The structure of the paper, summarized below, emphasizes the key ideas and
results brought to light in this work. In the rest of this section we present
the notation we use and we give some basic preliminary facts that are used
several times in the sequel. The second section deals with the identification
of a weak boundedness property of a set with respect to a closed and convex
cone and we call this dually cone-boundedness since its very definition
supposes that the images of the set through all the positive functionals of
the cone are bounded from below (as subsets of the set of real numbers). In
this section, we also compare the new concept with the usual notions of
cone-boundedness (understood either in the norm topology or in the weak
topology), as well as with other, seemingly unrelated notions, such as
hyperbolicity and pseudo-hyperbolicity of a set. Several examples distinguish
these notions from one another. The second section shows that several conic
cancellation rules that were recently proposed in literature are available
under the weaker assumption of dually cone-boundedness. Moreover, it is shown
that, in a certain sense, this type of conic boundedness is the weakest one
that still ensures the validity of the conclusion of the main cancellation
rules under scrutiny. The fourth section is dedicated to applications and is
divided into two subsections. The first subsection presents some statements
concerning the Pareto minimality under topological properties discussed
before. The second subsection, based on the second section, considers an
embedding procedure in the sense of R\aa dstr\"{o}m (\cite{Rad}) for dually
cone-bounded sets. We investigate and compare the R\aa dstr\"{o}m spaces
associated to three classes of sets and we focus on the possibility to
introduce a meaningful topology on the R\aa dstr\"{o}m space generated by the
dually cone-bounded sets. In particular, we study the closedness of the
ordering cone in this space as well as the nonemptiness of its interior. Some
differences with respect to the case of cone-bounded sets are emphasized, as
well. This allows us to conclude that several recent statements related to the
role of cone-boundedness in some classes of optimization problems can be
extended to the case of dually cone-boundedness. The paper ends with a short
concluding section.

\bigskip

Next we present the main notations and the objects we are using in this work.
Let $\left(  X,\left\Vert \cdot\right\Vert \right)  $ be a normed vector space
and $\left(  X^{\ast},\left\Vert \cdot\right\Vert _{\ast}\right)  $ its
topological dual. For a subset of $X$ or $X^{\ast}$ the standard symbols
$\operatorname{cl},$ $\operatorname*{int},$ $\operatorname*{conv},$
$\operatorname*{cone}$ stand for the closure, interior (with respect to the
norm, i.e., in the strong topology), the convex hull, and the conic hull,
respectively. The closure in the the weak$^{\ast}$ topology of $X^{\ast}$ is
denoted by $\operatorname{cl}_{w^{\ast}}$ (and similarly for other
topologies), while the bracket $\left\langle \cdot,\cdot\right\rangle $
denotes the duality mapping between $X$ and $X^{\ast}.$ By $D_{X}$ we
designate the closed unit ball (the discus) of $X.$ If $A,B$ are nonempty
subsets of $X$, the excess from $A$ to $B$ is defined by $e(A,B)=\sup_{x\in
A}d(x,B),$ where, as customary, $d(x,B)=\inf\left\{  \left\Vert x-b\right\Vert
\mid b\in B\right\}  .$ The Hausdorff-Pompeiu pseudo-distance between $A$ and
$B$ is $h\left(  A,B\right)  =\max\left\{  e(A,B),e(B,A)\right\}  .$

For convex cones $Q\subset X$ and $E\subset X^{\ast}$ we denote their positive
polars as
\[
Q^{+}=\left\{  x^{\ast}\in X^{\ast}\mid x^{\ast}\left(  x\right)  \geq0,\text{
}\forall x\in Q\right\}
\]
and
\[
E^{+}=\left\{  x\in X\mid x^{\ast}\left(  x\right)  \geq0,\text{ }\forall
x^{\ast}\in E\right\}  .
\]
We put $Q^{-}$ and $E^{-}$ for $-Q^{+}$ and $-E^{+},$ respectively. By the
Bipolar Theorem (see \cite[Theorem 1.1.9]{CZ}), $Q^{++}=\operatorname{cl}Q$
and $E^{++}=\operatorname{cl}_{w^{\ast}}E.$ Recall as well that $Q$ is said to
be pointed if $Q\cap-Q=\left\{  0\right\}  .$

\bigskip

In what follows, we consider $K\subset X$ as closed and convex cone. For a
nonempty subset $A$ of $X$ we denote $\tilde{A}=\operatorname{cl}%
\operatorname*{conv}\left(  A+K\right)  .$ Using repeatedly the elementary
relations
\[
\operatorname{cl}\left(  M+N\right)  =\operatorname{cl}\left(
\operatorname{cl}M+N\right)  \,\text{and }\operatorname*{conv}\left(
M+N\right)  =\operatorname*{conv}M+\operatorname*{conv}N,
\]
available for all nonempty sets $M,N\subset X,$ observe that $\tilde
{A}+K=\tilde{A}$ and $\widetilde{A+B}=\operatorname*{cl}\left(  \tilde
{A}+\tilde{B}\right)  $ for all $\emptyset\neq A,B\subset X.$ Indeed, the
former equality is proved by the following chain of inclusions:
\begin{align*}
\tilde{A}+K  &  \supset\tilde{A}=\operatorname*{cl}\operatorname*{conv}\left(
A+K\right)  =\operatorname{cl}\left(  \operatorname*{conv}A+K\right)
=\operatorname{cl}\left(  \operatorname*{conv}\left(  A+K\right)  +K\right) \\
&  =\operatorname{cl}\left(  \operatorname{cl}\operatorname*{conv}\left(
A+K\right)  +K\right)  \supset\tilde{A}+K.
\end{align*}
The latter relation is justified as follows:
\begin{align*}
\widetilde{A+B}  &  =\operatorname{cl}\operatorname*{conv}\left(
A+B+K\right)  =\operatorname{cl}\left(  \operatorname*{conv}\left(
A+K\right)  +\operatorname*{conv}\left(  B+K\right)  \right) \\
&  =\operatorname{cl}\left(  \operatorname{cl}\operatorname*{conv}\left(
A+K\right)  +\operatorname{cl}\operatorname*{conv}\left(  B+K\right)  \right)
=\operatorname{cl}\left(  \tilde{A}+\tilde{B}\right)  .
\end{align*}

\section{Dually cone-boundedness}

Consider, as mentioned before, that $(X,\Vert\cdot\Vert)$ is a normed vector
space, and let $K\subset X$ be a closed and convex cone. We are interested in
the study of those sets $\emptyset\neq A\subset X$ that satisfy the given
property that for every $x^{\ast}\in K^{+}$ the set $x^{\ast}(A)$ is bounded
from below. We will show that, on the one hand, this notion is strictly weaker
than the standard notions of $K-$boundedness with respect to the norm or weak
topology (see \cite[Definition 3.1, p. 13]{Luc}), and, on the other hand, it
allows us to extend several meaningful results with potential applications to
various types of optimization problems.

Firstly, we compare the newly proposed property with the $K-$boundedness with
respect to the strong topology and to this aim we begin with a technical
lemma, which we prove here for the sake of completeness.

\begin{lm}
\label{lm_dist_con}For all $x\in X$, we have
\[
d(x,K)=\sup\{-x^{\ast}(x)\mid x^{\ast}\in K^{+}\cap D_{X^{\ast}}\}.
\]

\end{lm}

\noindent\textbf{Proof. }If $d(x,K)=0$, the claim is trivial. Take $x^{\ast
}\in K^{+}$ with $\Vert x^{\ast}\Vert\leq1,$ and $k\in K$. Then, since
$x^{\ast}(k)\geq0$,
\[
-x^{\ast}(x)\leq-x^{\ast}(x-k)\leq\Vert x^{\ast}\Vert\cdot\Vert x-k\Vert
\leq\Vert x-k\Vert.
\]
Taking the infimum over $k\in K$ gives $-x^{\ast}(x)\leq d(x,K).$ Now taking
the supremum over all such $x^{\ast}$ yields one\ inequality.

Denote $\rho=d(x,K)>0$ and fix $\varepsilon\in(0,\rho)$. Then $K\cap D\left(
x,\rho-\varepsilon\right)  =\emptyset$. Since both sets are convex and,
moreover, the interior of the second one is nonempty, by the classical
separation theorem, there exist $x^{\ast}\in X^{\ast}\setminus\{0\}$ and
$\alpha\in\mathbb{R}$ such that
\[
\sup_{u\in D\left(  x,\rho-\varepsilon\right)  }x^{\ast}(u)\leq x^{\ast
}(k),\text{ }\forall k\in K.
\]
By a standard argument, $x^{\ast}\in K^{+}$. Then,
\[
x^{\ast}(x)+(\rho-\varepsilon)\Vert x^{\ast}\Vert=\sup_{\Vert u\Vert\leq
\rho-\varepsilon}x^{\ast}(x+u)=\sup_{u\in D\left(  x,\rho-\varepsilon\right)
}x^{\ast}(u)\leq0,
\]
whence
\[
-x^{\ast}(x)\geq(\rho-\varepsilon)\Vert x^{\ast}\Vert.
\]
Consider $\overline{x}^{\ast}=\left\Vert x^{\ast}\right\Vert ^{-1}x^{\ast}\in
K^{+}$. Then $-\overline{x}^{\ast}(x)\geq\rho-\varepsilon,$ whence
\[
\sup\{-x^{\ast}(x)\mid x^{\ast}\in K^{+}\cap D_{X^{\ast}}\}\geq\rho
-\varepsilon.
\]
Letting $\varepsilon\rightarrow0$, we obtain the other
needed\ inequality.\hfill$\square$

\begin{pr}
\label{prop1}Let $(X,\Vert\cdot\Vert)$ be a normed vector space, and let
$K\subset X$ be a closed and convex cone. Let $A\subset X$ be a nonempty set.
Consider the following assertions:

(i) there exists $\ell\geq0$ such that
\[
A\subset\ell D_{X}+K;
\]

(ii) there exists $\rho\in\mathbb{R}$ such that for every $x^{\ast}\in
K^{+}\cap D_{X^{\ast}},$ $\inf x^{\ast}\left(  A\right)  \geq\rho;$

(iii) for every $x^{\ast}\in K^{+}$ the set $x^{\ast}(A)$ is bounded from below.

Then (i)$\Leftrightarrow$(ii)$\Rightarrow$(iii).
\end{pr}

\noindent\textbf{Proof. }The fact that (i) implies (ii) with $\rho=-\ell$ is
clear. In order to prove that (ii) implies (i)$,$ we use the identity from
Lemma \ref{lm_dist_con}. For any $a\in A$ we have
\[
d(a,K)=\sup\{-x^{\ast}(a)\mid x^{\ast}\in K^{+}\cap D_{X^{\ast}}\}\leq
\sup_{x^{\ast}\in K^{+}\cap D_{X^{\ast}}}\sup\left\{  -x^{\ast}(u)\mid u\in
A\right\}  )\leq-\rho.
\]
Setting $\ell=-\rho\geq0$, we get $d(a,K)\leq\ell$ for all $a\in A$, i.e.
\[
A\subset\left(  \ell+\delta\right)  D_{X}+K,\text{ }\forall\delta>0,
\]
and the conclusion follows.

The implication (ii)$\Rightarrow$(iii) is obvious.\hfill$\square$

\bigskip

A set $A$ that satisfies property (i) from Proposition \ref{prop1} is called
in literature cone-bounded with respect to $K,$ or $K-$bounded (with respect
to the strong topology: see \cite{Luc}). A key point of our investigation is
that, in general, (iii) does not imply (i) and (ii) from Proposition
\ref{prop1}. In order to see that, we have the following example.

\begin{examp}
\label{example_b_b}Let $\left(  X,\left\Vert \cdot\right\Vert \right)
=\left(  c_{0},\left\Vert \cdot\right\Vert _{\infty}\right)  $ and $K=\left(
c_{0}\right)  _{+}\ $(the nonnegative cone). Then $\left(  X^{\ast},\left\Vert
\cdot\right\Vert _{\ast}\right)  =\left(  \ell^{1},\left\Vert \cdot\right\Vert
_{1}\right)  $ and $K^{+}=\ell_{+}^{1}.$

For every $k\in\mathbb{N}\setminus\left\{  0\right\}  $, define $a^{(k)}\in
c_{0}$ by
\[
a_{n}^{(k)}=%
\begin{cases}
k^{2}, & 1\leq n\leq k,\\
-k, & n=k+1,\\
0, & n\geq k+2,
\end{cases}
\]
and let $A=\{a^{(k)}\mid k\in\mathbb{N}\setminus\left\{  0\right\}  \}$.

Take $x^{\ast}=(x_{n})\in K^{+}$ with $\left\Vert x^{\ast}\right\Vert =1$. We
show that the set $x^{\ast}(A)$ is bounded below, which means that (iii) holds
in this particular case.

Indeed,
\[
x^{\ast}(a^{(k)})=\sum_{n=1}^{k}k^{2}x_{n}-kx_{k+1}=k^{2}\sum_{n=1}^{k}%
x_{n}-kx_{k+1}.
\]
Since $\sum_{n=1}^{\infty}x_{n}=1$ and all $x_{n}$ are nonnegative, there
exists $k_{x^{\ast}}$ such that for all $k\geq k_{x^{\ast}}$, $x_{k+1}%
\leq2^{-1}$ and $\sum_{n=1}^{k}x_{n}\geq2^{-1}.$ Then for $k\geq k_{x^{\ast}}%
$,
\[
x^{\ast}(a^{(k)})\geq2^{-1}\left(  k^{2}-k\right)  >0.
\]
Thus $x^{\ast}(a^{(k)})$ can be negative only for finitely many $k,$ whence
$\inf x^{\ast}(a^{(k)})>-\infty$.

Now we show that (i) does not hold. For $a^{(k)}$, the only negative
coordinate is $a_{k+1}^{(k)}=-k$, hence, for all $k,$ $d(a^{(k)},K)=k,$ and
consequently there exists no $\ell>0$ with $A\subset\ell D_{X}+K$.
\end{examp}

Another example, even in finite dimension, will be presented a little bit
later, after the description of a situation when the implication
(iii)$\Rightarrow$(i) in Proposition \ref{prop1} holds.

\begin{pr}
\label{prop2}If $X$ is finite dimensional normed vector space and $K$ is
finitely generated, then all three assertions from Proposition \ref{prop1} are equivalent.
\end{pr}

\noindent\textbf{Proof. }We have to prove only that (iii)$\Rightarrow$(i).
Since $X$ is finite dimensional, we identify $X^{\ast}$ with $X.$ Moreover,
since $K$ is finitely generated, $K^{+}$ is also finitely generated (see
\cite{Bes}). Let $u_{1},...,u_{n}\in X$ (where $n\in\mathbb{N}\setminus
\left\{  0\right\}  $) such that
\[
K^{+}=\operatorname*{cone}\left\{  u_{1},u_{2},...,u_{n}\right\}  .
\]
Then
\[
K=\left\{  a\in X\mid\left\langle a,u_{i}\right\rangle \geq0,\text{ }\forall
i\in\overline{1,n}\right\}  =\left\{  a\in X\mid\left\langle a,-u_{i}%
\right\rangle \leq0,\text{ }\forall i\in\overline{1,n}\right\}  .
\]
Assumption (ii) is equivalent to $\inf u_{i}\left(  A\right)  >-\infty$ for
all $i\in\overline{1,n},$ which is to say $\sup\left(  -u_{i}\left(  A\right)
\right)  <\infty$ for all $i\in\overline{1,n}.$ Denote by $c$ the value
$\max\left\{  0,\left(  \sup\left(  -u_{i}\left(  A\right)  \right)  \right)
_{i\in\overline{1,n}}\right\}  $ and observe that
\[
A\subset\left\{  x\in X\mid-u_{i}\left(  x\right)  \leq c,\text{ }\forall
i\in\overline{1,n}\right\}  .
\]
Now, by the Hoffman's Lemma (see \cite{Hof}), there is a constant $\rho>0$
such that for all $x\in A,$ $d\left(  x,K\right)  \leq\rho c.$ Consequently,
$(i)$ holds.\hfill$\square$

\begin{cor}
If $X$ is a normed vector space of dimension $2,$ then all three assertions
from Proposition \ref{prop1} are equivalent.
\end{cor}

\noindent\textbf{Proof. }In dimension $2$ all the closed convex cones are
finitely generated.\hfill$\square$

\begin{rmk}
Of course, under the assumption that $K^{+}$ is finitely generated, the space
$X$ could be infinite dimensional and Proposition \ref{prop2} is still true
(with basically the same proof). However, this requirement on $K^{+}$ is
rather strong, since it implies that $K$ is defined by a finite number of
linear inequalities, making the result essentially a finite dimensional one.
\end{rmk}

In order to show that one cannot avoid the requirement on the cone to be
finitely generated in Proposition \ref{prop2} we present the following example
which is to be found in \cite{Gos}.

\begin{examp}
\label{ex_Gos}Consider $X=\mathbb{R}^{3}$ and let $K$ be the ice-cream cone,
that is,
\[
K=\left\{  \left(  x,y,z\right)  \in\mathbb{R}^{3}\mid z\geq\sqrt{x^{2}+y^{2}%
}\right\}  .
\]
Take
\[
A=\left\{  x_{n}\mid x_{n}=\frac{n}{1-\cos\frac{1}{n}}\left(  \cos\frac{1}%
{n},\sin\frac{1}{n},\cos\frac{1}{n}\right)  ,\forall n\in\mathbb{N}%
\setminus\left\{  0\right\}  \right\}  .
\]
Then $K$ is not finitely generated and $K^{+}=K.$ However, for all $u\in
K^{+},$ $\left\langle u,x_{n}\right\rangle \geq0$ for all $n$ large enough,
while $d\left(  x_{n},K\right)  =\sqrt{2^{-1}}n$ for all $n\in\mathbb{N}%
\setminus\left\{  0\right\}  .$
\end{examp}

Examples \ref{example_b_b} and \ref{ex_Gos} allow us to define the following
distinct (strictly weaker) notion than $K-$boundedness.

\begin{df}
A set $A$ that satisfies property $(iii)$ from Proposition \ref{prop1} is
called dually cone-bounded with respect to $K$ or dually $K-$bounded.
\end{df}

\begin{rmk}
Observe that $A$ is dually $K-$bounded if and only if $x^{\ast}\left(
A\right)  \subset\mathbb{R}$ is $[0,\infty)-$bounded (that is, bounded from
below) for all $x^{\ast}\in K^{+}.$
\end{rmk}

Next we establish several results that allow a precise classification of the
dually cone-boundedness among several concepts of boundedness with respect to
a cone.

\begin{pr}
\label{prop_xs_weak}Let $\emptyset\neq A\subset X.$ The following assertions
are equivalent:

(i) $A$ is dually $K-$bounded;

(ii) for all $x^{\ast}\in X^{\ast}$, $e\left(  x^{\ast}\left(  A\right)
,x^{\ast}\left(  K\right)  \right)  \in\mathbb{R}.$
\end{pr}

\noindent\textbf{Proof. }Assume first that $A$ is dually $K-$bounded and take
$x^{\ast}\in X^{\ast}.$ We distinguish several cases. If $x^{\ast}\in K^{+},$
then $x^{\ast}\left(  K\right)  $ is either $\left\{  0\right\}  $ or
$[0,\infty),$ so for all $a\in A,$ we have that either $d\left(  x^{\ast
}\left(  a\right)  ,x^{\ast}\left(  K\right)  \right)  =0$ (if $x^{\ast
}\left(  a\right)  \geq0$)$,$ or $\left\vert x^{\ast}\left(  a\right)
\right\vert \leq\left\vert \inf x^{\ast}\left(  A\right)  \right\vert $ (if
$x^{\ast}\left(  a\right)  <0$)$.$ So, $\sup_{a\in A}d\left(  x^{\ast}\left(
a\right)  ,x^{\ast}\left(  K\right)  \right)  \leq\left\vert \inf x^{\ast
}\left(  A\right)  \right\vert \in\mathbb{R}.$ If $x^{\ast}\in K^{-}$ (whence
$-x^{\ast}\in K^{+}$) then $x^{\ast}\left(  K\right)  $ is either $\left\{
0\right\}  $ or $(-\infty,0],$ so for all $a\in A,$ either $d\left(  x^{\ast
}\left(  a\right)  ,x^{\ast}\left(  K\right)  \right)  =0$ (if $x^{\ast
}\left(  a\right)  \leq0$)$,$ or $x^{\ast}\left(  a\right)  =-\left(
-x^{\ast}\right)  \left(  a\right)  =\left\vert \left(  -x^{\ast}\right)
\left(  a\right)  \right\vert \leq\left\vert \inf\left(  -x^{\ast}\right)
\left(  A\right)  \right\vert $ (if $x^{\ast}\left(  a\right)  >0$). So, again
$\sup_{a\in A}d\left(  x^{\ast}\left(  a\right)  ,x^{\ast}\left(  K\right)
\right)  \leq\left\vert \inf\left(  -x^{\ast}\right)  \left(  A\right)
\right\vert \in\mathbb{R}.$ Finally, if $x^{\ast}\notin K^{+}\cup K^{-},$ then
$x^{\ast}\left(  K\right)  =\mathbb{R}$ and, consequently, $\sup_{a\in
A}d\left(  x^{\ast}\left(  a\right)  ,x^{\ast}\left(  K\right)  \right)  =0.$
So, (ii) is true.

Assume now that (ii) holds and take $x^{\ast}\in K^{+}.$ Let $\alpha>0$ such
that $\sup_{a\in A}d\left(  x^{\ast}\left(  a\right)  ,x^{\ast}\left(
K\right)  \right)  <\alpha.$ Then for all $a\in A,$ there is $k_{a}\in K$ such
that $\left\vert x^{\ast}\left(  a\right)  -x^{\ast}\left(  k_{a}\right)
\right\vert <\alpha.$ We deduce that
\[
-\alpha<x^{\ast}\left(  a\right)  -x^{\ast}\left(  k_{a}\right)  \leq x^{\ast
}\left(  a\right)  ,\text{ }\forall a\in A,
\]
whence $\inf x^{\ast}\left(  A\right)  \in\mathbb{R},$ meaning that (i) is
true.\hfill$\square$

\begin{rmk}
Proposition \ref{prop_xs_weak} and Example \ref{ex_Gos} mean that dually
$K-$boundedness is strictly weaker the $K-$boundedness with respect to the
weak topology, which according to \cite[p. 13]{Luc} means that for all weak
neighborhood $U$ of $0$ one can find $\alpha_{U}>0$ such that $A\subset
\alpha_{U}U+K.$ Of course, taking into account the form of the weak
neighborhoods of $0,$ this means that for all $n\in\mathbb{N}\setminus\left\{
0\right\}  $ and $\left(  x_{i}^{\ast}\right)  _{i\in\overline{1,n}}\subset
X^{\ast}$ there is $\alpha>0$ (depending on $n$ and $\left(  x_{i}^{\ast
}\right)  _{i\in\overline{1,n}}$) such that for all $a\in A$ there is
$k_{a}\in K$ with $\left\vert x_{i}^{\ast}\left(  a\right)  -x_{i}^{\ast
}\left(  k_{a}\right)  \right\vert <\alpha$ for all $i\in\overline{1,n}.$
\end{rmk}

Define $m_{A}:X^{\ast}\rightarrow\mathbb{R\cup}\left\{  -\infty\right\}  $
given by $m_{A}\left(  x^{\ast}\right)  =\inf x^{\ast}\left(  A\right)  .$
Assumption (iii) from Proposition \ref{prop1} reads as $m_{A}\left(
K^{+}\right)  \subset\mathbb{R}.$ Clearly, $m_{A}$ is a concave $w^{\ast}%
-$upper semicontinuous function. We further analyze the properties of this
function related to our issues by considering another way to look at the two
distinct properties of boundedness with respect to a cone as given by
Proposition \ref{prop1}. To this aim, recall that for a nonempty closed convex
set $C\subset X$, the recession cone can be defined as
\[
\operatorname*{rec}C=%
{\displaystyle\bigcap\limits_{\alpha>0}}
\alpha\left(  C-c\right)
\]
for every $c\in C$ (see, for instance, \cite[p. 6]{CZ}). The barrier cone of
$C$ is defined as
\[
\operatorname*{bar}C=\left\{  x^{\ast}\in X^{\ast}\mid\sup_{c\in C}x^{\ast
}\left(  c\right)  \in\mathbb{R}\right\}  .
\]
It is well known that $\operatorname{cl}_{w^{\ast}}\operatorname*{bar}%
C=\left(  \operatorname*{rec}C\right)  ^{-},$ and $\left(  \operatorname*{bar}%
C\right)  ^{-}=\operatorname*{rec}C$ (see \cite{CZ}).

After recalling the notation $\tilde{A}=\operatorname{cl}\operatorname*{conv}%
\left(  A+K\right)  ,$ some remarks are in order.

\begin{rmk}
\label{rmk1}Observe that $A\subset X$ is (dually) cone-bounded with respect to
$K$ if and only if the set $\tilde{A}$ is (dually) cone-bounded with respect
to $K$. Moreover, for all $x^{\ast}\in K^{+},$ $\inf x^{\ast}\left(  A\right)
=\inf x^{\ast}\left(  \tilde{A}\right)  $ (i.e., $m_{A}=m_{\tilde{A}}$ on
$K^{+}$) and if $A$ is dually cone-bounded then $\operatorname*{rec}\tilde
{A}=K.$ While the former claim is quite obvious, let us prove the latter one.
Indeed, if $k\in K,$ $a\in\tilde{A}$ and $t>0,$ then $a+tk\in\tilde
{A}+K=\tilde{A},$ which means that $K\subset\operatorname*{rec}\tilde{A}.$ Let
consider $u\in\operatorname*{rec}\tilde{A}.$ Then for an (all) element(s)
$a\in A$ and $t>0,$ $a+tu\in\tilde{A},$ so for every $x^{\ast}\in K^{+},$
$x^{\ast}\left(  a\right)  +tx^{\ast}\left(  u\right)  \geq\inf x^{\ast
}\left(  \tilde{A}\right)  \in\mathbb{R},$ which forces $x^{\ast}\left(
u\right)  \geq0.$ Because this is true for all $x^{\ast}\in K^{+},$ we get
that $u\in K,$ whence $\operatorname*{rec}\tilde{A}\subset K.$
\end{rmk}

\begin{rmk}
\label{rmk2}From the previous remark, we conclude that if $A$ is dually
cone-bounded, then $A\subset X$ is cone-bounded if and only if there is
$\ell>0$ such that $\tilde{A}\subset\ell D_{X}+\operatorname*{rec}\tilde{A}.$
A set having the latter property is said to be a hyperbolic set (see
\cite{Gos} and the references therein).
\end{rmk}

\begin{rmk}
\label{rmk3}Using \cite[Proposition 1.4]{Gos} we see that $\tilde{A}$ is
hyperbolic if and only if $m_{\tilde{A}}$ is bounded from below on every
bounded subset of $\left(  \operatorname*{rec}\tilde{A}\right)  ^{-}$ and if
and only if the restriction of $m_{\tilde{A}}$ to $\left(  \operatorname*{rec}%
\tilde{A}\right)  ^{-}$ is strongly continuous (that is, with respect to the
norm topology) at $0.$ In view of Remark \ref{rmk1}, we see that these
statements are in relation with the equivalence (i)$\Leftrightarrow$(ii) from
Proposition \ref{prop1}.
\end{rmk}

\begin{rmk}
Using similar types of arguments as in the above chain of remarks, we can draw
as well a link between Proposition \ref{prop2} and \cite[Proposition
4]{Bair85}. However, the setting, context, motivation and proofs of our
results are different from those quoted from \cite{Gos} and \cite{Bair85}.
\end{rmk}

Using the objects discussed in Remarks \ref{rmk1} and \ref{rmk2} we get the
following result.

\begin{pr}
\label{prop_b_h}A set $A\subset X$ is $K-$bounded if and only if
$\operatorname*{rec}\tilde{A}=K$ and $\tilde{A}$ is hyperbolic$.$
\end{pr}

\noindent\textbf{Proof. }Assume that $A$ is $K-$bounded. In particular, it is
dually $K-$bounded. By Remark \ref{rmk1}, we get that $\operatorname*{rec}%
\tilde{A}=K.$ Again, by the fact that $\tilde{A}$ is $K-$bounded, there is
$\ell>0$ such that $\tilde{A}\subset\ell D_{X}+K=\ell D_{X}%
+\operatorname*{rec}\tilde{A},$ whence $\tilde{A}$ is hyperbolic. The converse
is obvious.\hfill$\square$

\bigskip

Further characterizations of dually $K-$boundedness are to be observed, as follows.

\begin{pr}
A set $A\subset X$ is dually $K-$bounded if and only if $\operatorname*{bar}%
\tilde{A}=K^{-}.$
\end{pr}

\noindent\textbf{Proof. }If $A\subset X$ is dually $K-$bounded, then it is
easy to see that $K^{-}\subset\operatorname*{bar}\tilde{A}.$ Let $x^{\ast
}\notin K^{-}.$ By a standard separation result, there is $k\in K$ such that
$x^{\ast}\left(  k\right)  >0.$ We deduce that for all $\lambda>0,$ $x^{\ast
}\left(  \lambda k\right)  =\lambda x^{\ast}\left(  k\right)  >0$ and
$\sup_{x\in\tilde{A}}x^{\ast}\left(  x\right)  =\infty.$ The equality
$\operatorname*{bar}\tilde{A}=K^{-}$ now follows. The converse is
easy.\hfill$\square$

\bigskip

Using \cite[Definitions (p. 606), Proposition 1.2]{Gos}, we say, in our
context, that a convex set $C\subset X$ is pseudo-hyperbolic if
$\operatorname*{bar}C=\left(  \operatorname*{rec}C\right)  ^{-}.$ We get the
next consequence.

\begin{cor}
\label{cor_db_ph}A set $A\subset X$ is dually $K-$bounded if and only if
$\operatorname*{rec}\tilde{A}=K$ and $\tilde{A}$ is pseudo-hyperbolic.
\end{cor}

\begin{rmk}
At this point, we would like to underline an issue concerning the statement
that the norm-closedness of the barrier cone of a convex set (whence its
pseudo-hyperbolicity) implies hyperbolicity of that set. Although a proof of
this implication appeared in Bair's 1983 paper \cite{Bair83}, Bair later
acknowledged in his 1985 paper \cite{Bair85} that the result holds only in
dimension two. Moreover, Goossens in \cite[p. 607]{Gos} provided a
counterexample showing that the implication fails in dimension three (this is,
in fact, Example \ref{ex_Gos} quoted above). The distinction between these two
notions in the general setting is also evident from our preceding results
Proposition \ref{prop_b_h} and Corollary \ref{cor_db_ph}. We mention that the
relationships between these concepts and other related ones are used and
further explored in several recent papers, among which we cite here
\cite{Zaf}, \cite{EE}, and \cite{DL}.
\end{rmk}

Based on the above results and Proposition \ref{prop2} we get several equivalences.

\begin{pr}
Suppose that $X$ is finite dimensional and $K$ is finite generated. Then the
following assertions are equivalent:

(i) $A$ is $K-$bounded;

(ii) $A$ is dually $K-$bounded;

(iii) $\operatorname*{rec}\tilde{A}=K$ and $\tilde{A}$ is hyperbolic;

(iv) $\operatorname*{rec}\tilde{A}=K$ and $\tilde{A}$ is pseudo-hyperbolic.
\end{pr}

As already mentioned, $m_{A}$ is a concave $w^{\ast}-$upper semicontinuous
function. If it would be also $w^{\ast}-$lower semicontinuous (whence
$w^{\ast}-$continuous) on $K^{+}$, then by the fact that $K^{+}\cap
D_{X^{\ast}}$ is $w^{\ast}-$compact, (ii) from Proposition \ref{prop1} would
be verified. This gives us the motivation to record two chains of implications
firstly in finite dimension and then in infinite dimension. The notation is
that from above.

\begin{cor}
\label{cor_fin}Let $(X,\Vert\cdot\Vert)$ be a finite dimensional normed vector
space and $\emptyset\neq A\subset X$. Each of the following assertions implies
the following one:

(i) there exists a bounded set $B\subset X$ such that $A=B+K;$

(ii) $m_{A}$ is Lipschitz on $K^{+};$

(iii) $m_{A}$ is continuous on $K^{+};$

(iv) there is $\rho\in\mathbb{R}$ such that for every $x^{\ast}\in K^{+},$
$\inf x^{\ast}\left(  A\right)  \geq\rho\left\Vert x^{\ast}\right\Vert ;$

(v) there exists $\ell\geq0$ such that $A\subset\ell D_{X}+K;$

(vi) $m_{A}$ has finite values on $K^{+}.$

Moreover, (iv) and (v) are equivalent.
\end{cor}

\noindent\textbf{Proof. }(i)$\Rightarrow$(ii) is the first part of \cite[Lemma
4.13]{DFS}; (ii)$\Rightarrow$(iii) is obvious; (iii)$\Rightarrow$(iv) is true
by the application of the Weierstrass Theorem on the compact set $D_{X^{\ast}%
}\cap K^{+};$ (iv)$\Rightarrow$(v) is (i)$\Leftrightarrow$(ii) from
Proposition \ref{prop1}; (v)$\Rightarrow$(vi) is straightforward.\hfill
$\square$

\begin{cor}
\label{cor_inf}Let $(X,\Vert\cdot\Vert)$ be Banach space and $\emptyset\neq
A\subset X$. Consider the following assertions:

(i) there exists a compact set $B\subset X$ such that $A=B+K;$

(ii) $m_{A}$ is $w^{\ast}-$continuous on $K^{+};$

(iii) $m_{A}$ is $w^{\ast}-$sequentially continuous on $K^{+};$

(iv) there is $\rho>0$ such that for every $x^{\ast}\in K^{+},$ $\inf x^{\ast
}\left(  A\right)  \geq\rho\left\Vert x^{\ast}\right\Vert ;$

(v) there exists $\ell\geq0$ such that $A\subset\ell D_{X}+K;$

(vi) $m_{A}$ has finite values on $K^{+}.$

Then (i)$\Rightarrow$(iii), (ii)$\Rightarrow$(iv)$\Leftrightarrow
$(v)$\Rightarrow$(vi).
\end{cor}

\noindent\textbf{Proof. }(i)$\Rightarrow$(iii) is the second part of
\cite[Lemma 4.13]{DFS} and this is the only place where the completeness of
$X$ is needed; (ii)$\Rightarrow$(iv) is true by the application of the
Weierstrass Theorem on the $w^{\ast}-$compact set $D_{X^{\ast}}\cap K^{+};$
(iv)$\Leftrightarrow$(v) is (i)$\Leftrightarrow$(ii) from Proposition
\ref{prop1}; (v)$\Rightarrow$(vi) is straightforward.\hfill$\square$

\begin{rmk}
Notice that the assumptions (i) from the above two corollaries are
generalizations of the so-called Motzkin decomposition: see, for instance,
\cite{GLT}.
\end{rmk}

Some further examples illustrate the above results. The next example shows
that the strong continuity of $m_{A}$ (and therefore its $w^{\ast}-$lower
semicontinuity) is not implied by (iii) from Corollary \ref{cor_inf}.

\begin{examp}
Let $\left(  X,\left\Vert \cdot\right\Vert \right)  =\left(  c_{0},\left\Vert
\cdot\right\Vert _{\infty}\right)  $ and $K=\left(  c_{0}\right)  _{+}\ $(the
nonnegative cone). Consequently, $\left(  X^{\ast},\left\Vert \cdot\right\Vert
_{\ast}\right)  =\left(  \ell^{1},\left\Vert \cdot\right\Vert _{1}\right)  $
and $K^{+}=\ell_{+}^{1}.$ Let
\[
A=\{a^{(n)}=2^{n}e_{n}\mid n\in\mathbb{N}\setminus\left\{  0\right\}
\}\subset K,
\]
where $e_{n}$ is the $n$-th unit vector in $c_{0}$. Clearly, $A$ is $K-$bounded.

Define $x^{\ast}=(2^{-k})_{k}\in\ell_{+}^{1}$ and for all $n\in\mathbb{N}%
\setminus\left\{  0\right\}  ,$ $x_{n}^{\ast}\in\ell_{+}^{1}$ given by
\[
(x_{n}^{\ast})_{k}=%
\begin{cases}
0, & k=n,\\
2^{-k}, & k\neq n.
\end{cases}
\]
Then
\[
\Vert x_{n}^{\ast}-x^{\ast}\Vert_{1}=\sum_{k=1}^{\infty}|(x_{n}^{\ast}%
)_{k}-x_{k}|=|0-2^{-n}|=2^{-n}\rightarrow0,
\]
so $x_{n}^{\ast}\rightarrow x^{\ast}$ strongly.

Then for all $z^{\ast}=\left(  z_{k}\right)  _{k}\in\ell^{1},$
\[
\inf z^{\ast}(A)=\inf_{n}z^{\ast}(a^{(n)})=\inf_{n}\left(  2^{n}z_{n}\right)
.
\]
So, $\inf x^{\ast}(A)=1$ and $\inf x_{n}^{\ast}(A)=0$ for all $n.$ Thus,
although $x_{n}^{\ast}\rightarrow x^{\ast}$ strongly in $K^{+}$,
\[
\inf x_{n}^{\ast}(A)=0\not \rightarrow 1=\inf x^{\ast}(A).
\]
So, $m_{A}$ is not continuous in the norm topology, whence it is not $w^{\ast
}-$continuous.
\end{examp}

We present now an example, on $\mathbb{R}^{2},$ where $m_{A}$ is continuous
but not Lipschitz on $K^{+},$ whence (iii) does not imply (ii) in Corollary
\ref{cor_fin}.

\begin{cor}
Let $X=\mathbb{R}^{2}$ with the Euclidean norm and $K=\mathbb{R}_{+}^{2}$ (the
positive orthant), so $K^{+}=K$. Take
\[
A=\{(t,\,t^{-1})\mid t>0\}\subset K.
\]
In particular, $A$ is $K-$bounded with $\ell=0.$ For $x^{\ast}=(a,b)\in K^{+}%
$,
\[
m_{A}\left(  x^{\ast}\right)  =\inf_{t>0}\left(  at+bt^{-1}\right)  ,
\]
which gives
\[
m_{A}\left(  x^{\ast}\right)  =\left\{
\begin{array}
[c]{l}%
2\sqrt{ab},\text{if }a,b>0,\\
0,\text{if }a=0\ \text{or }b=0.
\end{array}
\right.
\]
So $m_{A}$ is continuous on $K^{+}$. However, $m_{A}$ is not Lipschitz on
$K^{+}$. Indeed, take $x^{\ast}=(0,\,1)$ and $y_{n}^{\ast}=(n^{-1},\,1)\ $for
all $n\in\mathbb{N}\setminus\left\{  0\right\}  .$ Then%
\[
\left\Vert x^{\ast}-y_{n}^{\ast}\right\Vert =n^{-1}\text{ and }\left\vert
m_{A}(x^{\ast})-m_{A}(y_{n}^{\ast})\right\vert =2\sqrt{n^{-1}},\text{ }\forall
n>0.
\]
Thus
\[
\frac{\left\vert m_{A}(x^{\ast})-m_{A}(y_{n}^{\ast})\right\vert }{\left\Vert
x^{\ast}-y_{n}^{\ast}\right\Vert }=2\sqrt{n}\rightarrow\infty,
\]
which means that there is no global Lipschitz constant for $m_{A}$ on $K^{+}.$
\end{cor}

At this point, it is worth considering a closely related side discussion about
the representation of convex cones with nonempty interior. Let $U\subset
X^{\ast}$ be a $w^{\ast}-$compact set and define
\begin{equation}
C=\left\{  x\in X\mid x^{\ast}\left(  x\right)  >0,\text{ }\forall x^{\ast}\in
U\right\}  . \label{the_set_C}%
\end{equation}
The question is to decide whether $C$ is open, an issue which was studied very
recently in \cite{LP}. We consider $m_{U}^{\ast}:X\rightarrow\mathbb{R}$ given
by $m_{U}^{\ast}\left(  x\right)  =\inf_{x^{\ast}\in U}x^{\ast}\left(
x\right)  ,$ that is, the same kind of function as $m_{A}$ discussed before,
but this time defined on the primal space. Taking into account the assumption
of $U$ we see that $m_{U}^{\ast}$ is well-defined and $C=\left\{  x\in X\mid
m_{U}^{\ast}\left(  x\right)  >0\right\}  .$ Clearly, if $m_{U}^{\ast}$ is
continuous, then $C$ is open. This is true, for instance, if $X$ is a Banach
space since in such a situation $U$ is strongly bounded, whence for all
$x,y\in X$ and $x^{\ast}\in U$ (denoting by $M$ the boundedness constant of
$U$) we have%
\[
x^{\ast}\left(  x\right)  =x^{\ast}\left(  y\right)  +x^{\ast}\left(
x-y\right)  \leq x^{\ast}\left(  y\right)  +\left\Vert x^{\ast}\right\Vert
\cdot\left\Vert x-y\right\Vert \leq x^{\ast}\left(  y\right)  +M\left\Vert
x-y\right\Vert ,
\]
whence
\[
m_{U}^{\ast}\left(  x\right)  \leq m_{U}^{\ast}\left(  y\right)  +M\left\Vert
x-y\right\Vert ,
\]
showing that $m_{U}^{\ast}$ is $M-$Lipschitz. However, in general (without the
completeness of $X$) the openness of $C$ affirmed by "if part" in
\cite[Proposition 3.12]{LP} is no longer true, as the next example proves.

\begin{examp}
Take $\left(  X,\left\Vert \cdot\right\Vert \right)  =\left(  c_{00}%
,\left\Vert \cdot\right\Vert _{\infty}\right)  $ for which $\left(  X^{\ast
},\left\Vert \cdot\right\Vert _{\ast}\right)  =\left(  \ell^{1},\left\Vert
\cdot\right\Vert _{1}\right)  .$ Consider for all $n\in\mathbb{N}%
\setminus\left\{  0\right\}  $ the element $x_{n}^{\ast}=\left(
0,...,0,n,0,...\right)  \in\ell^{1}$, where $n$ is in the $n-$th position.
Denote by $e_{1}^{\ast}=\left(  1,0,...,0,...\right)  \in\ell^{1}.$ The set
\[
U=\left\{  e_{1}^{\ast}\right\}  \cup\left\{  e_{1}^{\ast}+x_{n}^{\ast}\mid
n\in\mathbb{N}\setminus\left\{  0\right\}  \right\}
\]
is $w^{\ast}-$compact, a fact which is easy to see if one takes into account
that $\left(  x_{n}^{\ast}\right)  $ is $w^{\ast}-$convergent towards $0.$
Consider the set $C$ defined by (\ref{the_set_C}). Clearly, $x=\left(
1,0,...,0,...\right)  \in C.$ Take, for all $n\in\mathbb{N}\setminus\left\{
0\right\}  ,$ $x_{n}=\left(  1,0,...,-\sqrt{n^{-1}},0,...\right)  \in X,$
where $-\sqrt{n^{-1}}$ is in the $n-$th position. We can see that $\left(
e_{1}^{\ast}+x_{n}^{\ast}\right)  \left(  x_{n}\right)  =1-\sqrt{n}<0$ for all
$n\geq2.$ So, $x_{n}\notin C$ for $n\geq2.$ But, $\left(  x_{n}\right)  $
converges towards $x$ in $X,$ hence $x$ does not belong to the interior of
$C.$
\end{examp}

\section{Cancellation rules for dually cone-bounded sets}

In this section, we highlight an important property enjoyed by dually
cone-boundedness. More specifically, firstly we employ this concept in getting
some new conic versions of the celebrated R\aa dstr\"{o}m cancellation lemma
(see \cite[Lemma 1]{Rad} and \cite[Lemma 2.1]{Sch}). We mention that
generalizations of this result and some of its further versions in a conic
framework were obtained in \cite{DFcancellation}. Secondly, we show that the
property characterizing this class of sets is, in a sense, the minimal
requirement for obtaining results of this type.

\begin{pr}
\label{pr_Rad_conic_dual}Suppose that $A,B,C\subset X\ $are nonempty sets such
that $C$ is dually $K-$bounded and
\[
A+C\subset\operatorname{cl}\left(  C+B+K\right)  .
\]
Then
\[
A\subset\operatorname{cl}\operatorname*{conv}\left(  B+K\right)  .
\]

\end{pr}

\noindent\textbf{Proof. }The assumption lets us write
\[
A+C\subset C+B+K+D\left(  0,\varepsilon\right)  ,\text{ }\forall
\varepsilon>0.
\]
Let $x^{\ast}\in K^{+}\setminus\left\{  0\right\}  .$ Then $x^{\ast}\left(
A+C\right)  \subset x^{\ast}\left(  C+B+K+D\left(  0,\varepsilon\right)
\right)  ,$ which means that
\[
x^{\ast}\left(  A\right)  +x^{\ast}\left(  C\right)  \subset x^{\ast}\left(
C\right)  +x^{\ast}\left(  B\right)  +x^{\ast}\left(  D\left(  0,\varepsilon
\right)  \right)  +[0,\infty),\text{ }\forall\varepsilon>0.
\]
This relation allows to write%
\[
\inf x^{\ast}\left(  A\right)  +\inf x^{\ast}\left(  C\right)  \geq\inf
x^{\ast}\left(  C\right)  +\inf x^{\ast}\left(  B\right)  -\varepsilon
\left\Vert x^{\ast}\right\Vert ,\text{ }\forall\varepsilon>0,
\]
where, as usual, one applies the convention $-\infty+a=-\infty$ for all
$a\in\lbrack-\infty,\infty)$.

Since $C$ is dually $K-$bounded, $\inf x^{\ast}\left(  C\right)  \in
\mathbb{R},$ so, letting $\varepsilon\rightarrow0,$ we get $\inf x^{\ast
}\left(  A\right)  \geq\inf x^{\ast}\left(  B\right)  .$

Assume now, by way of contradiction, that there is $\overline{x}\in
A\setminus\operatorname{cl}\operatorname*{conv}\left(  B+K\right)  .$ Then we
strongly separate the point $\overline{x}$ from the set $\operatorname{cl}%
\operatorname*{conv}\left(  B+K\right)  ,$ so we get $x^{\ast}\in X^{\ast
}\setminus\left\{  0\right\}  $ and $\alpha\in\mathbb{R}$ such that%
\[
x^{\ast}\left(  \overline{x}\right)  <\alpha<x^{\ast}\left(  b+k\right)
,\text{ }\forall b\in B,\text{ }\forall k\in K.
\]
By a standard argument, $x^{\ast}\in K^{+}\setminus\left\{  0\right\}  $ and
\[
x^{\ast}\left(  \overline{x}\right)  <\alpha<x^{\ast}\left(  b\right)  ,\text{
}\forall b\in B.
\]
But this implies as well that $\inf x^{\ast}\left(  A\right)  \leq x^{\ast
}\left(  \overline{x}\right)  <\alpha\leq\inf x^{\ast}\left(  B\right)  ,$
contradicting the above step of the proof.\hfill$\square$

\bigskip

In view of Examples \ref{example_b_b} and \ref{ex_Gos}, this result
generalizes \cite[Proposition 2.2]{DFcancellation}. Another conic cancellation
law, under different assumptions, reads as follows.

\begin{pr}
\label{prop_canc_sol}Suppose that $A,B,C\subset X\ $are nonempty sets and
$\operatorname*{int}K\neq\emptyset$. If $A$ and $C$ are dually $K-$bounded
and
\[
A+C\subset C+B+\operatorname*{int}K,
\]
then
\[
A\subset\operatorname*{conv}B+\operatorname*{int}K.
\]

\end{pr}

\noindent\textbf{Proof. }As above, for all $x^{\ast}\in K^{+}\setminus\left\{
0\right\}  ,$
\[
x^{\ast}\left(  A\right)  +x^{\ast}\left(  C\right)  \subset x^{\ast}\left(
C\right)  +x^{\ast}\left(  B\right)  +(0,\infty).
\]

Since $A\,$and $C$ are dually $K-$bounded, $\inf x^{\ast}\left(  A\right)
+\inf x^{\ast}\left(  C\right)  \in\mathbb{R}$ and the above relation assures
that
\[
\inf x^{\ast}\left(  A\right)  +\inf x^{\ast}\left(  C\right)  >\inf x^{\ast
}\left(  C\right)  +\inf x^{\ast}\left(  B\right)  ,
\]
whence $\inf x^{\ast}\left(  A\right)  >\inf x^{\ast}\left(  B\right)  .$

Assume now, by way of contradiction, that there is $\overline{x}\in
A\setminus\left(  \operatorname*{conv}B+\operatorname*{int}K\right)  .$ Then,
again by separation, we get $x^{\ast}\in X^{\ast}\setminus\left\{  0\right\}
$ and $\alpha\in\mathbb{R}$ such that%
\[
x^{\ast}\left(  \overline{x}\right)  <x^{\ast}\left(  b+k\right)  ,\text{
}\forall b\in B,\text{ }\forall k\in\operatorname*{int}K.
\]
Surely, $x^{\ast}\in K^{+}\setminus\left\{  0\right\}  $ and
\[
x^{\ast}\left(  \overline{x}\right)  \leq x^{\ast}\left(  b\right)  ,\text{
}\forall b\in B.
\]
This means as well that $\inf x^{\ast}\left(  A\right)  \leq x^{\ast}\left(
\overline{x}\right)  \leq\inf x^{\ast}\left(  B\right)  .$ This is a
contradiction, and the conclusion follows.\hfill$\square$

\bigskip

When a comparison is made between the above two results and the corresponding
\cite[Propositions 2.2 and 2.6]{DFcancellation}, it can be observed that the
notion of dually cone-boundedness provides a very natural framework for conic
cancellation laws. Indeed, one can remark that the concept of weakly
$K-$compactness defined in \cite{DFcancellation} implies that of dually
$K-$boundedness (see \cite[Lemma 2.5]{DFcancellation}).

Next, we show that the conic cancellation rules cannot be further extended
beyond the class of dually cone-bounded sets.

\begin{lm}
\label{lm_equi_star}Let $A,B\subset X$ be nonempty sets, and $E\subset
X^{\ast}$ be a $w^{\ast}-$closed and convex cone. Then the following
assertions are equivalent:

(i) $\inf x^{\ast}\left(  A\right)  \geq\inf x^{\ast}\left(  B\right)  $ for
all $x^{\ast}\in E;$

(ii) $A\subset\operatorname{cl}\operatorname*{conv}\left(  B+E^{+}\right)  .$
\end{lm}

\noindent\textbf{Proof. }Assume (i) as hypothesis. If (ii) were false, then we
would find $a\in A\setminus\operatorname{cl}\operatorname*{conv}\left(
B+E^{+}\right)  .$ Then there would exist $x^{\ast}\in X^{\ast}$ such that
$x^{\ast}\left(  a\right)  <\inf x^{\ast}\left(  B\right)  +\inf x^{\ast
}\left(  E^{+}\right)  .$ By a well-known argument, we get that $x^{\ast
}\left(  m\right)  \geq0,$ for all $m\in E^{+},$ whence $x^{\ast}\in
E^{++}=\operatorname{cl}_{w^{\ast}}E=E.$ Therefore $\inf x^{\ast}\left(
E^{+}\right)  =0$ and the inequality $\inf x^{\ast}\left(  A\right)  \leq
x^{\ast}\left(  a\right)  <\inf x^{\ast}\left(  B\right)  $ contradicts the assumption.

Suppose now that (ii) is true. Then for all $x^{\ast}\in E,$%
\[
\inf x^{\ast}\left(  A\right)  \geq\inf x^{\ast}\left(  B\right)  +\inf
x^{\ast}\left(  E^{+}\right)  \geq\inf x^{\ast}\left(  B\right)  ,
\]
so (i) is true as well.\hfill$\square$

\begin{pr}
\label{pr_equiv_ph2}Assume that $C\subset X$ is a set such that $\tilde{C}$ is
pseudo-hyperbolic. Let $A,B\subset X.$ Then the following assertions are equivalent:

(i) $A+C\subset\operatorname{cl}\operatorname*{conv}\left(  B+\tilde
{C}\right)  ;$

(ii) $A\subset\operatorname{cl}\operatorname*{conv}\left(
B+\operatorname*{rec}\tilde{C}\right)  .$
\end{pr}

\noindent\textbf{Proof. }We proceed by equivalence. Since $\widetilde
{M+N}=\operatorname{cl}\left(  \tilde{M}+\tilde{N}\right)  $ for any sets
$M,N\subset X,$ the relation $A+C\subset\operatorname{cl}\operatorname*{conv}%
\left(  B+\tilde{C}\right)  $ is equivalent to $A+\tilde{C}\subset
\operatorname{cl}\operatorname*{conv}\left(  B+\tilde{C}\right)  .$ The latter
inclusion is equivalent (again via a separation argument) to
\[
\inf x^{\ast}\left(  A\right)  +\inf x^{\ast}\left(  \tilde{C}\right)
\geq\inf x^{\ast}\left(  B\right)  +\inf x^{\ast}\left(  \tilde{C}\right)
,\text{ }\forall x^{\ast}\in X^{\ast},
\]
and then further equivalent to
\[
\inf x^{\ast}\left(  A\right)  \geq\inf x^{\ast}\left(  B\right)  ,\text{
}\forall x^{\ast}\in-\operatorname*{bar}\tilde{C}.
\]
By Lemma \ref{lm_equi_star} (notice that $\operatorname*{bar}\tilde{C}$ is
$w^{\ast}-$closed because $\tilde{C}$ is pseudo-hyperbolic) the above relation
is equivalent to
\[
A\subset\operatorname{cl}\operatorname*{conv}\left(  B+\left(
-\operatorname*{bar}\tilde{C}\right)  ^{+}\right)  =\operatorname{cl}%
\operatorname*{conv}\left(  B+\operatorname*{rec}\tilde{C}\right)  .
\]
The proof is complete.\hfill$\square$

\begin{pr}
Assume that $C\subset X$ is a set such that $\tilde{C}$ is pseudo-hyperbolic.
Then the following assertions are equivalent:

(i) For all $\emptyset\neq A,B\subset X$ the inclusion
\[
A+\tilde{C}\subset\operatorname{cl}\operatorname*{conv}\left(  B+\tilde
{C}+K\right)
\]
implies that
\[
A\subset\operatorname{cl}\operatorname*{conv}\left(  B+K\right)  ;
\]

(ii) $C$ is dually $K-$bounded.
\end{pr}

\noindent\textbf{Proof. }We prove that (i)$\Rightarrow$(ii). The inclusion
\[
A+\tilde{C}\subset\operatorname{cl}\operatorname*{conv}\left(  B+\tilde
{C}+K\right)
\]
is the same with
\[
A+\tilde{C}\subset\operatorname{cl}\operatorname*{conv}\left(  B+\tilde
{C}\right)  ,
\]
and Proposition \ref{pr_equiv_ph2} shows that this is equivalent with
\[
A\subset\operatorname{cl}\operatorname*{conv}\left(  B+\operatorname*{rec}%
\tilde{C}\right)  .
\]
Taking here $A=\operatorname*{rec}\tilde{C}$ and $B=\left\{  0\right\}  $ the
above relation is true, whence so is the equivalent relation before it.
Applying the assumption in this case gives $\operatorname*{rec}\tilde
{C}\subset K.$ As observed in Remark \ref{rmk1}, the converse inclusion is
always true, whence $\operatorname*{rec}\tilde{C}=K.$ Taking into account that
$\tilde{C}$ is pseudo-hyperbolic and using Corollary \ref{cor_db_ph} we obtain
that $C$ is dually $K-$bounded.

The implication (ii)$\Rightarrow$(i) follows by the fact that dually
$K-$boundedness of $C$ is equivalent to dually $K-$boundedness of $\tilde{C}$
and from Proposition \ref{pr_Rad_conic_dual}.\hfill$\square$

\section{Applications}

In this section, we present several applications in which the notion of dually
cone-boundedness proves useful. In particular, we show that in various
optimization settings, one can replace stronger boundedness assumptions with
the central concept introduced in this work, while still preserving the
essential features of the corresponding results. We organize this section into
two subsections.

\subsection{Application to Pareto points}

First, we provide a link between the dually cone-boundedness and the Pareto
minimality. Let $A\subset X.$ Recall (see \cite{GRTZ}) that $a\in A$ is a
(weak) Pareto minimum point of $A$ with respect to $K$ if $\left(  A-a\right)
\cap-K=\left\{  0\right\}  $ (respectively, $\left(  A-a\right)
\cap-\operatorname*{int}K=\emptyset$ provided $\operatorname*{int}%
K\neq\emptyset$).

\begin{df}
In the above notation we say that a dually $K-$bounded set $A$ has the dually
$K-$minimality property if for all $x^{\ast}\in K^{+},$ there exists $\min
x^{\ast}\left(  A\right)  .$
\end{df}

\begin{rmk}
\label{rem_w_b}Clearly, if $x^{\ast}\left(  A\right)  $ is closed for all
$x^{\ast}\in K^{+},$ then $A$ is dually $K-$bounded if and only if it has the
dually $K-$minimality property.
\end{rmk}

We briefly recall now some concepts of generalized compactness. According to
\cite{Luc} a set $A\subset X$ is called $K-$compact (or compact with respect
to the cone $K$) if from any cover of $A$ with the sets of the form $U+K$,
where $U$ is open, one can extract a finite subcover of it. In
\cite{DF-coneComp} the following concept was introduced and studied: the set
$A$ is called $K-$sequentially compact if for any sequence $\left(
a_{n}\right)  \subset A$ there is a sequence $\left(  c_{n}\right)  \subset K$
such that the sequence $\left(  a_{n}-c_{n}\right)  $ has a convergent
subsequence towards an element of $A$. These two notions are also known under
the names of cone compactness and sequential cone compactness, respectively.
It was shown in \cite{DF-coneComp} that if $A$ is $K-$compact, then $A$ is
$K-$sequentially compact. For the converse, according to \cite[Theorem
2.11]{DF-coneComp}, if $A$ is $K-$sequentially compact and separable, then $A$
is $K-$compact, but without separability this implication fails (see
\cite{DFS}). Finally, again according to \cite{Luc}, $A$ is called $K-$convex
if $A+K$ is convex.

\begin{pr}
\label{prop_secv_b}Suppose that the set $A$ is $K-$sequentially compact and
$\operatorname*{int}K\neq\emptyset$. Then $A$ is dually $K-$bounded and $A$
has the dually $K-$minimality property.
\end{pr}

\noindent\textbf{Proof. }Let $x^{\ast}\in K^{+}\setminus\left\{  0\right\}  .$
Suppose, by way of contradiction, that there is $\left(  a_{n}\right)  \subset
A$ such that $x^{\ast}\left(  a_{n}\right)  \rightarrow-\infty.$ Then, since
$A$ is $K-$sequentially compact, we find a sequence $\left(  c_{n}\right)
\subset K$ such that $a_{n}-c_{n}\rightarrow a\in A.$ Consequently, $x^{\ast
}\left(  a_{n}\right)  =x^{\ast}\left(  a_{n}-c_{n}\right)  +x^{\ast}\left(
c_{n}\right)  \geq x^{\ast}\left(  a_{n}-c_{n}\right)  \rightarrow y^{\ast
}\left(  a\right)  ,$ which is a contradiction. So, $A$ is dually $K-$bounded.
Clearly, this also implies that $A+K$ is dually $K-$bounded.

We prove now that $x^{\ast}\left(  A+K\right)  $ is closed. Indeed, let
$\left(  a_{n}\right)  \subset A,$ $\left(  k_{n}\right)  \subset K$ such that
$x^{\ast}\left(  a_{n}+k_{n}\right)  \rightarrow\alpha\in\mathbb{R}.$ Again,
we find a sequence $\left(  c_{n}\right)  \subset K$ such that $a_{n}%
-c_{n}\rightarrow a\in A$, so $x^{\ast}\left(  a_{n}-c_{n}+c_{n}+k_{n}\right)
\rightarrow\alpha.$ This yields that $x^{\ast}\left(  c_{n}+k_{n}\right)
\rightarrow\alpha-x^{\ast}\left(  a\right)  .$ Since $\operatorname*{int}%
K\neq\emptyset,$ $x^{\ast}\left(  K\right)  =\mathbb{R}_{+},$ whence we deduce
that $\alpha\in x^{\ast}\left(  A\right)  +\mathbb{R}_{+}=x^{\ast}\left(
A+K\right)  ,$ and the thesis is proved. By the previous step of the proof, we
know that $A+K$ is dually $K-$bounded and following Remark \ref{rem_w_b} $A+K$
has the dually $K-$minimality property. But, obviously, this is equivalent to
say that $A$ has the dually $K-$minimality property. The conclusion
ensues.\hfill$\square$

\bigskip

The usefulness of these notions becomes clear in the next result.

\begin{pr}
Suppose that $\operatorname*{int}K\neq\emptyset.$ If there exists $x^{\ast}\in
K^{+}\setminus\left\{  0\right\}  $ such that there exists $\min x^{\ast
}\left(  A\right)  ,$ then $A$ has a weak Pareto minimum with respect to $K$.
If, moreover, $A$ is $K-$convex, then the converse holds.
\end{pr}

\noindent\textbf{Proof. }According to the assumption, there is $a\in A$ such
that for all $b\in A$, $x^{\ast}\left(  a\right)  \leq x^{\ast}\left(
b\right)  .$ If, there would exist $b\in A$ such that $b-a\in
-\operatorname*{int}K$ then $x^{\ast}\left(  b-a\right)  <0$ which is a contradiction.

Conversely, if $A$ is $K-$convex, and $a$ is a weak Pareto minimum of $A,$
then $\left(  a-\operatorname*{int}K\right)  \cap\operatorname{cl}\left(
A+K\right)  =\emptyset,$ so, by a separation theorem, there is $x^{\ast}\in
K^{+}\setminus\left\{  0\right\}  $ such that $x^{\ast}\left(  a\right)  \leq
x^{\ast}\left(  b\right)  $ for all $b\in A,$ so the conclusion of this
implication holds.\hfill$\square$

\subsection{Embedding results}

The cancellations rules derived in the preceding section lead to an embedding
result, allowing to see the family of dually cone-bounded sets as a cone in a
normed vector space. This provides a far-reaching extension of the original
R\aa dstr\"{o}m embedding theorem for convex compact sets in a normed vector
space from \cite{Rad}. Notice that such an extension was studied in
\cite{DSta} for the case of cone-bounded sets. Here we explain what parts and
what applications of the construction in \cite{DSta} could be transferred to
the more general case of dually cone-bounded sets.

The following sets are slight variations of those introduced in \cite{DSta}:%
\[
\mathcal{CB}_{K}\left(  X\right)  =\left\{  A\subset X\mid A\text{ is
nonempty, }K-\text{bounded}\right\}
\]
and
\[
\mathcal{C}_{K}\left(  X\right)  =\left\{  \tilde{A}\mid A\in\mathcal{CB}%
_{K}\left(  X\right)  \right\}  .
\]

Denote%
\[
\mathcal{CB}_{K^{+}}\left(  X\right)  =\left\{  A\subset X\mid A\text{ is
nonempty, dually }K-\text{bounded}\right\}
\]
and
\[
\mathcal{C}_{K^{+}}\left(  X\right)  =\left\{  \tilde{A}\mid A\in
\mathcal{CB}_{K^{+}}\left(  X\right)  \right\}  .
\]
Clearly $\mathcal{CB}_{K}\left(  X\right)  \subset\mathcal{CB}_{K^{+}}\left(
X\right)  $ and $\mathcal{C}_{K}\left(  X\right)  \subset\mathcal{C}_{K^{+}%
}\left(  X\right)  $ and both inclusion are strict in view of Goossens'
Example \ref{ex_Gos} and Example \ref{example_b_b}. These examples also
motivate us to introduce the set%
\[
\mathcal{C}_{G}\left(  X\right)  =\left\{  \tilde{A}\mid A\in\mathcal{CB}%
_{K^{+}}\left(  X\right)  \setminus\mathcal{CB}_{K}\left(  X\right)  \right\}
,
\]
where the subscript $G$ comes from Goossens.

For $A,B\subset X$ nonempty and dually $K-$bounded sets, define $\tilde
{A}\oplus\tilde{B}=\widetilde{A+B}.$ Remark that Proposition
\ref{pr_Rad_conic_dual} shows that if $A,B,C\in\mathcal{CB}_{K^{+}}\left(
X\right)  ,$ then the equality $\tilde{A}\oplus\tilde{C}=\tilde{C}\oplus
\tilde{B}$ implies $\tilde{A}=\tilde{B}.$ Consequently, $\left(
\mathcal{C}_{K^{+}}\left(  X\right)  ,\oplus\right)  $ is a commutative
monoid, with cancellation law, the neutral element being $K=\operatorname*{cl}%
\left(  \left\{  0\right\}  +K\right)  \in\mathcal{C}_{K^{+}}\left(  X\right)
.$

We define, as an external operation on $\mathcal{C}_{K^{+}}\left(  X\right)
,$ the multiplication with nonnegative scalars, denoted $\odot,$ as follows:%
\[
\lambda\odot\tilde{A}=\widetilde{\lambda A}=\left\{
\begin{array}
[c]{l}%
K,\text{ if }\lambda=0,\\
\lambda\tilde{A},\text{ if }\lambda>0.
\end{array}
\right.
\]
Moreover, for any $\tilde{A},\tilde{B}\in\mathcal{C}_{K^{+}}\left(  X\right)
$ and any $\lambda_{1},\lambda_{2}\geq0,$ one readily gets $\lambda_{1}%
\odot\left(  \tilde{A}\oplus\tilde{B}\right)  =\lambda_{1}\odot\tilde{A}%
\oplus\lambda_{1}\odot\tilde{B},$ $\left(  \lambda_{1}+\lambda_{2}\right)
\odot\tilde{A}=\lambda_{1}\odot\tilde{A}\oplus\lambda_{2}\odot\tilde{A},$
$\lambda_{1}\odot\left(  \lambda_{2}\odot\tilde{A}\right)  =\left(
\lambda_{1}\lambda_{2}\right)  \odot\tilde{A},$ $1\odot\tilde{A}=\tilde{A}.$

Now it is easy to see that all the assumptions from \cite[Theorem 1, A,
B]{Rad} are fulfilled and therefore $\mathcal{C}_{K^{+}}\left(  X\right)  $
can be embedded into a vector space which we denote $\mathcal{G}_{K^{+}%
}\left(  X\right)  .$ Following R\aa dstr\"{o}m's proof, this linear space
$\mathcal{G}_{K^{+}}\left(  X\right)  $ (the R\aa dstr\"{o}m space of
$\mathcal{C}_{K^{+}}\left(  X\right)  ,$ or associated to $\mathcal{C}_{K^{+}%
}\left(  X\right)  $) consists of all equivalence classes defined by the
equivalence relation $\sim$ on $\mathcal{C}_{K^{+}}\left(  X\right)
\times\mathcal{C}_{K^{+}}\left(  X\right)  $ given by
\[
\left(  \tilde{A},\tilde{B}\right)  \sim\left(  \tilde{C},\tilde{D}\right)
\iff\tilde{A}\oplus\tilde{D}=\tilde{B}\oplus\tilde{C}.
\]
More precisely, the addition on $\mathcal{G}_{K^{+}}\left(  X\right)  $ is
defined by
\[
\left\langle \tilde{A},\tilde{B}\right\rangle +\left\langle \tilde{C}%
,\tilde{D}\right\rangle =\left\langle \tilde{A}\oplus\tilde{C},\tilde{B}%
\oplus\tilde{D}\right\rangle ,
\]
while the multiplication by a scalar $\lambda\in\mathbb{R}$ is given by%
\[
\lambda\left\langle \tilde{A},\tilde{B}\right\rangle =\left\{
\begin{array}
[c]{l}%
\left\langle \lambda\odot\tilde{A},\lambda\odot\tilde{B}\right\rangle ,\text{
if }\lambda\geq0,\\
\left\langle -\lambda\odot\tilde{B},-\lambda\odot\tilde{A}\right\rangle
,\text{ if }\lambda<0,
\end{array}
\right.
\]
and these operations ensure a structure of a real linear space for
$\mathcal{G}_{K^{+}}\left(  X\right)  .$ The set $\mathcal{C}_{K^{+}}\left(
X\right)  $ is embedded into this vectorial space by the mapping
$\varphi:\mathcal{C}_{K^{+}}\left(  X\right)  \rightarrow\mathcal{G}_{K^{+}%
}\left(  X\right)  $ given by $\varphi\left(  \tilde{A}\right)  =\left\langle
\tilde{A},K\right\rangle ,$ and $\varphi\left(  \mathcal{C}_{K^{+}}\left(
X\right)  \right)  $ is a convex cone in $\mathcal{G}_{K^{+}}\left(  X\right)
.$

Furthermore, $\mathcal{G}_{K^{+}}\left(  X\right)  $ contains the convex
pointed cone defined as
\[
\mathcal{K}_{K^{+}}\left(  X\right)  =\left\{  \left\langle \tilde{A}%
,\tilde{B}\right\rangle \in\mathcal{G}_{K^{+}}\left(  X\right)  \mid\tilde
{A}\subset\tilde{B}\right\}  ,
\]
which we call the ordering cone of $\mathcal{G}_{K^{+}}\left(  X\right)  .$

Notice that the corresponding construction for $\mathcal{C}_{K}\left(
X\right)  $ was detailed in \cite{DSta} and the associated R\aa dstr\"{o}m
space was denoted by $\mathcal{G}_{K}\left(  X\right)  ,$ while its ordering
cone was designated by $\mathcal{K}_{K}\left(  X\right)  .$

\bigskip

Observe that if $\mathcal{C}_{G}\left(  X\right)  \neq\emptyset$ then
$\mathcal{C}_{G}\left(  X\right)  \cup\left\{  K\right\}  =\left(
\mathcal{C}_{K^{+}}\left(  X\right)  \setminus\mathcal{C}_{K}\left(  X\right)
\right)  \cup\left\{  K\right\}  $ is also a commutative monoid, with
cancellation law and the same neutral element, so in the same manner one can
construct its R\aa dstr\"{o}m space which we denote by $\mathcal{G}_{G}\left(
X\right)  .$ Some more remarks recorded in the following result are in order.

\begin{pr}
(i) $\mathcal{G}_{K}\left(  X\right)  $ isomorphic to a linear subspace of
$\mathcal{G}_{K^{+}}\left(  X\right)  $;

(ii) if $\mathcal{C}_{G}\left(  X\right)  \neq\emptyset$ then $\mathcal{G}%
_{G}\left(  X\right)  $ is isomorphic to $\mathcal{G}_{K^{+}}\left(  X\right)
.$
\end{pr}

\noindent\textbf{Proof. }(i) Define the linear application $\psi
:\mathcal{G}_{K}\left(  X\right)  \rightarrow\mathcal{G}_{K^{+}}\left(
X\right)  ,$%
\[
\psi\left(  \left\langle \tilde{A},\tilde{B}\right\rangle _{\mathcal{G}%
_{K}\left(  X\right)  }\right)  =\left\langle \tilde{A},\tilde{B}\right\rangle
_{\mathcal{G}_{K^{+}}\left(  X\right)  },
\]
where the subscripts indicate the space in which the class is considered. This
function is the isomorphism we are looking for.

(ii) Define the linear application $\omega:\mathcal{G}_{G}\left(  X\right)
\rightarrow\mathcal{G}_{K^{+}}\left(  X\right)  ,$%
\[
\omega\left(  \left\langle \tilde{A},\tilde{B}\right\rangle _{\mathcal{G}%
_{G}\left(  X\right)  }\right)  =\left\langle \tilde{A},\tilde{B}\right\rangle
_{\mathcal{G}_{K^{+}}\left(  X\right)  }.
\]
It is easy, indeed, to see that this function is linear and injective. If
$\left\langle \tilde{A},\tilde{B}\right\rangle \in\mathcal{G}_{K^{+}}\left(
X\right)  ,$ taking $\tilde{C}\in\mathcal{C}_{K^{+}}\left(  X\right)
\setminus\mathcal{C}_{K}\left(  X\right)  ,$ we have that $\tilde{A}%
\oplus\tilde{C},\tilde{B}\oplus\tilde{C}\in\mathcal{C}_{K^{+}}\left(
X\right)  \setminus\mathcal{C}_{K}\left(  X\right)  $,
\[
\left\langle \tilde{A}\oplus\tilde{C},\tilde{B}\oplus\tilde{C}\right\rangle
_{\mathcal{G}_{G}\left(  X\right)  }=\left\langle \tilde{A},\tilde
{B}\right\rangle _{\mathcal{G}_{G}\left(  X\right)  },
\]
and
\[
\omega\left(  \left\langle \tilde{A}\oplus\tilde{C},\tilde{B}\oplus\tilde
{C}\right\rangle _{\mathcal{G}_{G}\left(  X\right)  }\right)  =\left\langle
\tilde{A},\tilde{B}\right\rangle _{\mathcal{G}_{K^{+}}\left(  X\right)  }.
\]
Therefore the function is also surjective.\hfill$\square$

\begin{rmk}
Using the application $\psi$ defined in the proof of the above result, we also
have that $\psi\left(  \mathcal{K}_{K}\left(  X\right)  \right)
\subset\mathcal{K}_{K^{+}}\left(  X\right)  .$
\end{rmk}

In the case of $\mathcal{C}_{K}\left(  X\right)  $ we were able to perform the
extra step given by \cite[Theorem 1, C]{Rad}, that is to define a norm on
$\mathcal{G}_{K}\left(  X\right)  $ by $\left\Vert \left\langle \tilde
{A},\tilde{B}\right\rangle \right\Vert =h\left(  \tilde{A},\tilde{B}\right)
.$ Notice that for some $K-$bounded sets $\tilde{A},\tilde{B}$ the
Hausdorff-Pompeiu distance is finite and this norm actually extends the metric
$h$ from $\mathcal{C}_{K}\left(  X\right)  $ (see \cite[Lemma 3.1]{DSta}). In
the new case we are considering, that of $\mathcal{C}_{K^{+}}\left(  X\right)
,$ $h$ is not anymore with finite values. Notice that, actually, $h\left(
\tilde{A},K\right)  <\infty$ if and only if $A$ is $K-$bounded, whence
$h\left(  \tilde{A},K\right)  =\infty$ if $\tilde{A}\in\mathcal{C}_{K^{+}%
}\left(  X\right)  \setminus\mathcal{C}_{K}\left(  X\right)  .$ In this
respect we record the following result.

\begin{pr}
\label{pr_norm_inf}(i) If $\left\langle \tilde{A},\tilde{B}\right\rangle
=\left\langle \tilde{C},\tilde{D}\right\rangle \in\mathcal{G}_{K^{+}}\left(
X\right)  $ then $h\left(  \tilde{A},\tilde{B}\right)  =h\left(  \tilde
{C},\tilde{D}\right)  ;$

(ii) if $\tilde{A}\in\mathcal{C}_{K^{+}}\left(  X\right)  \setminus
\mathcal{C}_{K}\left(  X\right)  $ and $\tilde{B}\in\mathcal{C}_{K}\left(
X\right)  ,$ then $h\left(  \tilde{A},\tilde{B}\right)  =\infty.$
\end{pr}

\noindent\textbf{Proof. }(i)\textbf{ }It is enough to show that if $h\left(
\tilde{A},\tilde{B}\right)  <\infty,$ then $h\left(  \tilde{C},\tilde
{D}\right)  =h\left(  \tilde{A},\tilde{B}\right)  .$ Denote $h\left(
\tilde{A},\tilde{B}\right)  =\alpha\in\mathbb{R}.$ Then for all $\varepsilon
>0$
\[
\tilde{A}\subset\tilde{B}+\left(  \alpha+\varepsilon\right)  D_{X}\text{ and
}\tilde{B}\subset\tilde{A}+\left(  \alpha+\varepsilon\right)  D_{X}.
\]
Consequently,
\[
\tilde{A}+\tilde{D}\subset\widetilde{B+D}+\left(  \alpha+\varepsilon\right)
D_{X}\text{ and }\tilde{B}+\tilde{C}\subset\widetilde{A+C}+\left(
\alpha+\varepsilon\right)  D_{X},
\]
whence, by the fact that $\tilde{A}\oplus\tilde{D}=\tilde{B}\oplus\tilde{C},$%
\[
\tilde{B}+\tilde{C}\subset\operatorname{cl}\left(  \tilde{B}+\tilde{D}+\left(
\alpha+\varepsilon\right)  D_{X}\right)  \text{ and }\tilde{A}+\tilde
{D}\subset\operatorname{cl}\left(  \tilde{A}+\tilde{C}+\left(  \alpha
+\varepsilon\right)  D_{X}\right)  .
\]
By Proposition \ref{pr_Rad_conic_dual} we get%
\[
\tilde{C}\subset\operatorname{cl}\left(  \tilde{D}+\left(  \alpha
+\varepsilon\right)  D_{X}\right)  \subset\tilde{D}+\left(  \alpha
+2\varepsilon\right)  D_{X}\text{ and }\tilde{D}\subset\operatorname{cl}%
\left(  \tilde{C}+\left(  \alpha+\varepsilon\right)  D_{X}\right)
\subset\tilde{C}+\left(  \alpha+2\varepsilon\right)  D_{X}.
\]
We obtain that $h\left(  \tilde{C},\tilde{D}\right)  \leq\alpha+2\varepsilon$
and, letting $\varepsilon\rightarrow0,$ $h\left(  \tilde{C},\tilde{D}\right)
\leq\alpha=h\left(  \tilde{A},\tilde{B}\right)  .$ Changing the roles of
$\left\langle \tilde{A},\tilde{B}\right\rangle $ and $\left\langle \tilde
{C},\tilde{D}\right\rangle $ we deduce the reverse inequality, whence the equality.

(ii) If, by way of contradiction, $h\left(  \tilde{A},\tilde{B}\right)
<\infty,$ then for all $\varepsilon>0$%
\[
\tilde{A}\subset\tilde{B}+\left(  h\left(  \tilde{A},\tilde{B}\right)
+\varepsilon\right)  D_{X},
\]
and since the set in the right-hand side is $K-$bounded, we get that $A$ is
$K-$bounded, which is a contradiction.\hfill$\square$

\bigskip

Concerning the ways to introduce a topological structure on $\mathcal{G}%
_{K^{+}}\left(  X\right)  $ there are (at least) two possibilities.

One of them, which we briefly mention, is to consider the extended norm (that
is a norm with possible infinite values) of the form given above. Denote this
topology by $\tau_{\left\Vert \cdot\right\Vert }.$ We distinguish between two cases.

Firstly, if $\mathcal{C}_{G}\left(  X\right)  =\emptyset,$ then $\mathcal{G}%
_{K}\left(  X\right)  =\mathcal{G}_{K^{+}}\left(  X\right)  $ and, clearly,
$\tau_{\left\Vert \cdot\right\Vert }$ is a usual norm topology which was
studied in \cite{DSta} for $\mathcal{G}_{K}\left(  X\right)  .$ In this case
it was proven in the quoted work that $\operatorname*{int}\mathcal{K}%
_{K}\left(  X\right)  $ with respect to $\tau_{\left\Vert \cdot\right\Vert }$
is nonempty, provided $\operatorname*{int}K\neq\emptyset$ and, more precisely,%
\[
\operatorname*{int}\mathcal{K}_{K}\left(  X\right)  =\left\{  \left\langle
\tilde{A},\tilde{B}\right\rangle \in\mathcal{K}_{K}\left(  X\right)
\mid\exists\varepsilon>0\text{ s.t. }\tilde{A}+\varepsilon B_{X}\subset
\tilde{B}\right\}  .
\]

Secondly, if $\mathcal{C}_{G}\left(  X\right)  \neq\emptyset$ then, according
to Proposition \ref{pr_norm_inf} (ii), $\left\Vert \cdot\right\Vert $ assumes
infinite values, a situation thoroughly examined in \cite{Beer} for a general
setting. In this case, the interior of $\operatorname*{int}\mathcal{K}_{K^{+}%
}\left(  X\right)  $ is still nonempty with respect to $\tau_{\left\Vert
\cdot\right\Vert }$ if $\operatorname*{int}K\neq\emptyset,$ the proof being
similar to \cite[Proposition 3.2]{DSta}. Observe that in $\operatorname*{int}%
\mathcal{K}_{K^{+}}\left(  X\right)  $ there are elements of norm $\infty.$
Indeed, if one takes $C\in\mathcal{C}_{G}\left(  X\right)  $ and $B=C\cup
D_{X},$ then $B\in\mathcal{C}_{G}\left(  X\right)  $, and according to
Proposition \ref{pr_norm_inf} (ii) the element $\left\langle 2^{-1}%
\widetilde{D_{X}},\tilde{B}\right\rangle $ from $\operatorname*{int}%
\mathcal{K}_{K^{+}}\left(  X\right)  $ has the norm $\infty$. However, the
topology induced by such a norm is not a linear topology. For this reason, we
are motivated to find a proper topology on $\mathcal{G}_{K^{+}}\left(
X\right)  .$

\bigskip

Therefore, we consider the family of seminorms $\mathcal{P}=\left\{
p_{x^{\ast}}\mid x^{\ast}\in K^{+}\right\}  ,$ where $p_{x^{\ast}}%
:\mathcal{G}_{K^{+}}\left(  X\right)  \rightarrow\mathbb{R},$%
\[
p_{x^{\ast}}\left(  \left\langle \tilde{A},\tilde{B}\right\rangle \right)
=\left\vert \inf x^{\ast}\left(  \tilde{A}\right)  -\inf x^{\ast}\left(
\tilde{B}\right)  \right\vert =\left\vert \inf x^{\ast}\left(  A\right)  -\inf
x^{\ast}\left(  B\right)  \right\vert ,\text{ }\forall x^{\ast}\in K^{+},
\]
that endow $\mathcal{G}_{K^{+}}\left(  X\right)  $ with a locally convex
topology. We denote this topology by $\tau_{\mathcal{P}}.$ Of course, this is
a coarser topology than $\tau_{\left\Vert \cdot\right\Vert }.$

\bigskip

The main topological properties of $\mathcal{K}_{K^{+}}\left(  X\right)  $ in
$\tau_{\mathcal{P}}$ are given below.

\begin{pr}
\label{pr_con_cl}(i) The cone $\mathcal{K}_{K^{+}}\left(  X\right)  $ is
closed in $\tau_{\mathcal{P}}.$

(ii) If $K$ is pointed and $\mathcal{C}_{G}\left(  X\right)  \neq\emptyset$,
then the interior of $\mathcal{K}_{K^{+}}\left(  X\right)  $ in $\tau
_{\mathcal{P}}$ is empty.
\end{pr}

\noindent\textbf{Proof. }(i)\textbf{ }Let $\left\langle \tilde{A},\tilde
{B}\right\rangle \in\operatorname{cl}_{\tau_{\mathcal{P}}}\mathcal{K}_{K^{+}%
}\left(  X\right)  .$ Then, for all $x^{\ast}\in K^{+}$ and all $\varepsilon
>0,$ there is $\left\langle \tilde{C}_{x^{\ast},\varepsilon},\tilde
{D}_{x^{\ast},\varepsilon}\right\rangle \in\mathcal{G}_{K^{+}}\left(
X\right)  $ such that%
\[
\left\vert \inf x^{\ast}\left(  \tilde{C}_{x^{\ast},\varepsilon}\right)  -\inf
x^{\ast}\left(  \tilde{D}_{x^{\ast},\varepsilon}\right)  \right\vert
<\varepsilon
\]
and
\[
\left\langle \tilde{A},\tilde{B}\right\rangle +\left\langle \tilde{C}%
_{x^{\ast},\varepsilon},\tilde{D}_{x^{\ast},\varepsilon}\right\rangle
\in\mathcal{K}_{K^{+}}\left(  X\right)  .
\]
Taking into account the definition of $\mathcal{K}_{K^{+}}\left(  X\right)  $
we have the inclusion%
\[
\tilde{A}\oplus\tilde{C}_{x^{\ast},\varepsilon}\subset\tilde{B}\oplus\tilde
{D}_{x^{\ast},\varepsilon},
\]
whence%
\[
x^{\ast}\left(  \tilde{A}\right)  +x^{\ast}\left(  \tilde{C}_{x^{\ast
},\varepsilon}\right)  +\left(  0,\infty\right)  \subset x^{\ast}\left(
\tilde{B}\right)  +x^{\ast}\left(  \tilde{D}_{x^{\ast},\varepsilon}\right)
+\left(  -\varepsilon,\infty\right)  .
\]
Therefore, we can write
\begin{align*}
x^{\ast}\left(  \tilde{A}\right)  +\inf x^{\ast}\left(  \tilde{C}_{x^{\ast
},\varepsilon}\right)  +\left(  0,\infty\right)   &  \subset x^{\ast}\left(
\tilde{A}\right)  +x^{\ast}\left(  \tilde{C}_{x^{\ast},\varepsilon}\right)
+(-\varepsilon,0)+\left(  0,\infty\right) \\
&  \subset x^{\ast}\left(  \tilde{B}\right)  +x^{\ast}\left(  \tilde
{D}_{x^{\ast},\varepsilon}\right)  +(-\varepsilon,0)+\left(  -\varepsilon
,\infty\right) \\
&  \subset x^{\ast}\left(  \tilde{B}\right)  +\inf x^{\ast}\left(  \tilde
{D}_{x^{\ast},\varepsilon}\right)  +(-\varepsilon,0)+\left(  -\varepsilon
,\infty\right) \\
&  \subset x^{\ast}\left(  \tilde{B}\right)  +\inf x^{\ast}\left(  \tilde
{C}_{x^{\ast},\varepsilon}\right)  +\left(  -2\varepsilon,\varepsilon\right)
+\left(  -\varepsilon,\infty\right) \\
&  =x^{\ast}\left(  \tilde{B}\right)  +\inf x^{\ast}\left(  \tilde{C}%
_{x^{\ast},\varepsilon}\right)  +\left(  -3\varepsilon,\infty\right)  .
\end{align*}
So, for all $\varepsilon>0,$%
\[
x^{\ast}\left(  \tilde{A}\right)  \subset x^{\ast}\left(  \tilde{B}\right)
+\left(  -3\varepsilon,\infty\right)  ,
\]
and this gives%
\[
\inf x^{\ast}\left(  \tilde{A}\right)  \geq\inf x^{\ast}\left(  \tilde
{B}\right)  .
\]
Since the latter inequality holds for all $x^{\ast}\in K^{+},$ a simple
separation argument shows that $\tilde{A}\subset\tilde{B},$ so $\left\langle
\tilde{A},\tilde{B}\right\rangle \in\mathcal{K}_{K^{+}}\left(  X\right)  .$

(ii) Suppose, by way of contradiction, that the conclusion is not true. Then
there is $\left\langle \tilde{A},\tilde{B}\right\rangle $ in the interior of
$\mathcal{K}_{K^{+}}\left(  X\right)  $ with respect to $\tau_{\mathcal{P}}.$
Then one can find $\varepsilon>0,$ $n\in\mathbb{N}\setminus\left\{  0\right\}
,$ $\left(  x_{i}^{\ast}\right)  _{i\in\overline{1,n}}\subset K^{+}$ such that
for all $\left\langle \tilde{C},\tilde{D}\right\rangle \in\mathcal{G}_{K^{+}%
}\left(  X\right)  $ the relation
\[
\left\vert \inf x_{i}^{\ast}\left(  \tilde{C}\right)  -\inf x_{i}^{\ast
}\left(  \tilde{D}\right)  -\inf x_{i}^{\ast}\left(  \tilde{A}\right)  +\inf
x_{i}^{\ast}\left(  \tilde{B}\right)  \right\vert <\varepsilon,\text{ }\forall
i\in\overline{1,n}%
\]
implies that $\tilde{C}\subset\tilde{D}.$

Take now $\tilde{C}=\operatorname{cl}\left(  \tilde{A}+\bigcap_{i\in
\overline{1,n}}\ker x_{i}^{\ast}\right)  $ and $\tilde{D}=\tilde{B}.$ Then the
above relation is obviously true, whence
\[
\tilde{A}+\bigcap_{i\in\overline{1,n}}\ker x_{i}^{\ast}\subset\tilde{B}.
\]
Take $u\in\bigcap_{i\in\overline{1,n}}\ker x_{i}^{\ast}.$ Since $\tilde{A}%
\neq\emptyset,$ we get that $u\in\operatorname*{rec}\tilde{B}.$ But, according
to Remark \ref{rmk1}, $\operatorname*{rec}\tilde{B}=K.$ As $\bigcap
_{i\in\overline{1,n}}\ker x_{i}^{\ast}$ is a linear subspace, we infer that
\[
\bigcap_{i\in\overline{1,n}}\ker x_{i}^{\ast}\subset K\cap-K.
\]
As $K$ is pointed, we get that $\bigcap_{i\in\overline{1,n}}\ker x_{i}^{\ast
}=\left\{  0\right\}  .$ Therefore the application
\[
x\mapsto\left(  x_{1}^{\ast}\left(  x\right)  ,...,x_{n}^{\ast}\left(
x\right)  \right)  \in\mathbb{R}^{n}%
\]
is a linear injection having as range a finite dimensional space. This means
that $X$ is also finite dimensional.

Moreover, for all $x^{\ast}\in X^{\ast},$
\[
\left\{  0\right\}  =\bigcap_{i\in\overline{1,n}}\ker x_{i}^{\ast}\subset\ker
x^{\ast},
\]
and by a well-known result (see, e.g., \cite[Lemma 3.2]{Brezis}), $x^{\ast}$
can be written as a linear combination of $\left(  x_{i}^{\ast}\right)
_{i\in\overline{1,n}},$ so $X^{\ast}$ is generated, as linear space, by
$\left(  x_{i}^{\ast}\right)  _{i\in\overline{1,n}}.$ Denote
\[
K^{\symbol{94}}=\operatorname*{conv}\operatorname*{cone}\left\{  x_{i}^{\ast
}\mid i\in\overline{1,n}\right\}  ,
\]
which is a closed convex cone. The inclusion $K^{\symbol{94}}\subset K^{+}$ is
obvious. We show that the reverse inclusion is also true. Consider, for all
$\lambda>0,$ the set $K_{\lambda}=\left(  K^{\symbol{94}}\right)  ^{+}%
\cap\lambda D_{X}$ which is a dually $K-$bounded set (being actually bounded)
and take $A_{\lambda}=\tilde{A}\oplus\widetilde{K_{\lambda}}\in\mathcal{C}%
_{K^{+}}\left(  X\right)  $. It is easy to see that for all $i\in
\overline{1,n},$ $\inf x_{i}^{\ast}\left(  A_{\lambda}\right)  =\inf
x_{i}^{\ast}\left(  A\right)  ,$ so
\[
\left\vert \inf x_{i}^{\ast}\left(  A_{\lambda}\right)  -\inf x_{i}^{\ast
}\left(  \tilde{B}\right)  -\inf x_{i}^{\ast}\left(  \tilde{A}\right)  +\inf
x_{i}^{\ast}\left(  \tilde{B}\right)  \right\vert <\varepsilon,\text{ }\forall
i\in\overline{1,n}.
\]
We get that $A_{\lambda}\subset\tilde{B}$ for all $\lambda>0,$ whence
$\tilde{A}+\left(  K^{\symbol{94}}\right)  ^{+}\subset\tilde{B}.$ As above,
$\left(  K^{\symbol{94}}\right)  ^{+}\subset\operatorname*{rec}\tilde{B}=K,$
which gives $K^{+}\subset K^{\symbol{94}}.$ Therefore, finally,
$K^{\symbol{94}}=K^{+},$ an equality that shows that $K^{+}$ is finitely
generated. By Proposition \ref{prop2}, $\mathcal{C}_{G}\left(  X\right)
=\emptyset$, which is a contradiction with our assumptions. We conclude that
the interior of $\mathcal{K}_{K^{+}}\left(  X\right)  $ in $\tau_{\mathcal{P}%
}$ is empty.\hfill$\square$

\bigskip

At this point, in view of all these facts, we are able to affirm that all the
results and applications from the quoted work \cite{DSta} that take into
account only the algebraic structure of the embedding space and the closedness
of the ordering cone can be extended to the case of $\mathcal{G}_{K^{+}%
}\left(  X\right)  ,$ which implies that throughout these results dually
$K-$bounded sets may be used in place of $K-$bounded sets. Among these
extension there are as well some assertions concerning different types of
optimization problems.

\section{Concluding remarks}

In this work we devised a weak notion of boundedness of a subset of a normed
vector space with respect to a closed and convex cone. This notion is strictly
weaker than all the related notions found in the literature, yet it still has
enough strength to ensure the conclusion of some important results in
optimization. We also consider that this concept, which we have called dually
cone-boundedness, may prove useful in investigating other issues as well. For
instance, we think that it is compatible with the idea of extending some
concepts and results concerning the Pettis integrability of closed-valued
multifunctions (see, for instance, \cite{EH}). Another possible use of this
concept, which could be explored in future work, is in establishing weak
conditions for certain special types of domination properties in vector
optimization (see \cite{EE2}).

\bigskip

\noindent\textbf{Data availability.} This manuscript has no associated data.

\noindent\textbf{Disclosure statement. }The authors contributed equally to all
aspects of the study and to the preparation of the manuscript. No potential
conflict of interest was reported by the authors.


\bigskip

\begin{thebibliography}{99}                                                                                               %


\bibitem {Bair83}J. Bair, \emph{Liens entre le cone d'ouverture interne et
l'internat du cone asymptotique d'un convexe}, Bulletin de la Soci\'{e}t\'{e}
Math\'{e}matique de Belgique, S\'{e}rie B, 35 (1983), 177--187.

\bibitem {Bair85}J. Bair, \emph{Retour sur les caract\`{e}res ferm\'{e} et
relativement ouvert du c\^{o}ne-barri\`{e}re d'un convexe}, Commentationes
Mathematicae Universitatis Carolinae, 26 (1985), 315--318.

\bibitem {Bau}H.H. Bauschke, H. Ouyang, X. Wang, \emph{On angles between
convex cones}, Journal of Applied and Numerical Optimization, 4 (2022), 131--141.

\bibitem {Beer}G. Beer, \emph{Norms with infinite values}, Journal of Convex
Analysis, 22 (2015), 37--60.

\bibitem {Bes}D.P. Bertsekas, \emph{Convex Optimization Theory}, Athena
Scientific, Belmont, 2009.

\bibitem {Brezis}H. Brezis, \emph{Functional analysis, Sobolev spaces and
partial differential equations}, Springer, New York, 2011.

\bibitem {DL}D. D\"{o}rfler, A. L\"{o}hne, \emph{Convex sets approximable as
sum of a compact set and a cone}, Journal of Nonlinear and Variational
Analysis, 8 (2024), 681--689.

\bibitem {DF-coneComp}M. Durea, E.-A. Florea, \emph{Cone-compactness of a set
and applications to set-equilibrium problems}, Journal of Optimization Theory
and Applications, 200 (2024), 1286--1308.

\bibitem {DFcancellation}M. Durea, E.-A. Florea,\emph{ Conic cancellation laws
and some applications in set optimization}, Optimization, 74 (2025), 2657--2675.

\bibitem {DSta}M. Durea, E.-C. Stamate, \emph{An embedding result for a class
of epigraphical sets and applications to set optimization}, Journal of
Nonlinear and Variational Analysis, 9 (2025), 907--925.

\bibitem {DFS}M. Durea, E.-A. Florea, E.-C. Stamate, \emph{A common approach
to vector and set optimization under differentiability assumptions},
Optimization, DOI: 10.1080/02331934.2025.2597968.

\bibitem {EH}K. El Amri, C. Hess, \emph{On the Pettis integral of closed
valued multifunctions}, Set-Valued Analysis, 8 (2000), 329--360.

\bibitem {EE}E. Ernst, \emph{A converse of the Gale-Klee-Rockafellar Theorem:
continuity of convex functions at the boundary of their domains}, Proceedings
of the American Mathematical Society, 141 (2013), 3665--3672.

\bibitem {EE2}E. Ernst, \emph{Pareto efficiency for vector-valued functions
which are lower semi-continuous in the sense of Penot and Th\'{e}ra},
Optimization, DOI: 10.1080/02331934.2025.2589350.

\bibitem {Far}A. Farajzadeh, \emph{On maximal and minimal elements for sets
with respect to cones}, Optimization Letters, 18 (2024), 1533--1540.

\bibitem {GLT}M.A. Goberna, J.E. Mart\'{\i}nez-Legaz, M.I. Todorov, \emph{On
Motzkin decomposable sets and functions}, Journal of Mathematical Analysis and
Applications, 372 (2010), 525--537.

\bibitem {Gos}P. Goossens, \emph{Hyperbolic sets and asymptotes}, Journal of
Mathematical Analysis and Applications, 116 (1986), 604--618.

\bibitem {GRTZ}A. G\"{o}pfert, H. Riahi, C. Tammer, C. Z\u{a}linescu,
\emph{Variational Methods in Partially Ordered Spaces}, second edition,
Springer, Cham, 2023.

\bibitem {Hof}A.J. Hoffman, \emph{On approximate solutions of systems of
linear inequalities}, Journal of Research of the National Bureau of Standards,
49 (1952), 263--265.

\bibitem {LP}P. Leonetti, G. Principi, \emph{Representations of cones and
applications to decision theory}, Journal of Mathematical Analysis and
Applications, 550 (2025), paper no. 129561.

\bibitem {Luc}D.T. Luc, \emph{Theory of Vector Optimization}, Springer,
Berlin, 1989.

\bibitem {MP}J.E. Mart\'{\i}nez-Legaz, C. Pintea, \emph{Closed convex sets
with an open or closed Gauss range}, Mathematical Programming, 189 (2021), 433--454.

\bibitem {Rad}H. R\aa dstr\"{o}m, \emph{An embedding theorem for spaces of
convex sets}, Proceedings of the American Mathematical Society, 3 (1952), 165--169.

\bibitem {PT}J.-P. Penot, M. Th\'{e}ra, \emph{Semi-continuous mappings in
general topology}, Archiv der Mathematik, 38 (1982), 158--166.

\bibitem {Sch}K.D. Schmidt, \emph{Embedding theorems for classes of convex
sets}, Acta Applicandae Mathematicae, 5 (1986), 209--237.

\bibitem {Zaf}A. Zaffaroni, \emph{Convex radiant costarshaped sets and the
least sublinear gauge}, Journal of Convex Analysis, 20 (2013), 307--328.

\bibitem {CZ}C. Z\u{a}linescu, \emph{Convex Analysis in General Vector
Spaces}, World Scientific, Singapore, 2002.
\end{thebibliography}
\end{document}